\documentstyle[12pt,epsf,epsfig]{article}
\advance\hoffset by -7mm  
\makeatletter
\@addtoreset{equation}{section}
\makeatother
\renewcommand{\theequation}{\thesection.\arabic{equation}}
\renewcommand{\title}[1]{\null\vspace{25mm}

\noindent{\Large{\bf #1}}\vspace{10mm}

\noindent {\large By }}
\newcommand{\authors}[1]{\noindent{\large #1}\vspace{3mm}

}
\newcommand{\address}[1]{\noindent #1\vspace{5mm}

}
\renewcommand{\abstract}[1]{\vspace{19mm}

\noindent{\small{\em Abstract.} #1}\vspace{2mm}

} 
 \newcommand{\eg}{{\em e.g.$\,$}}
 \newcommand{\cf}{{\em cf.$\,$ }}
 \newcommand{\ie}{{\em i.e.$\,$}}
 \newcommand{\Iff}{{\em iff\/ $\,$}}
 
 \newcommand{\rhs}{right side }

 \newcommand{\QED}{\mbox{\rule[-1.5pt]{6pt}{10pt}}}
 
 \newcommand{\R}{I\!\!R}
 \newcommand{\C}{C\!\!\!\rule[.5pt]{.7pt}{6.5pt}\:\:}
 
 \newcommand{\N}{I\!\! N}
 \newcommand{\AAA}{{\cal A}}

 \newcommand{\DD}{{\cal D}}
 \newcommand{\GG}{{\cal G}}
 \newcommand{\HH}{{\cal H}}
 \newcommand{\KK}{{\cal K}}
 \newcommand{\JJ}{{\cal J}}
 
 \newcommand{\OO}{{\cal O}}
 \newcommand{\SS}{{\cal S}}
 \newcommand{\VV}{{\cal V}}
 \newcommand{\re}{{\rm Re\,}}
 \newcommand{\im}{{\rm Im\,}}
 \newcommand{\eps}{\varepsilon}
 \newtheorem{claim}{Claim}[section]
 \newtheorem{theorem}[claim]{Theorem}
\newtheorem{proposition}[claim]{Proposition}
\newtheorem{definition}[claim]{Definition}
\newtheorem{lemma}[claim]{Lemma}
\newtheorem{corollary}[claim]{Corollary}

\newtheorem{remark}[claim]{Remark}
\newtheorem{remarks}[claim]{Remarks}

\begin{document}
\title{Exponential Bounds on Curvature-Induced Resonances in a
Two-Dimensional Dirichlet Tube}
\authors{P. Duclos} 
\address{Centre de Physique Th\'eorique, CNRS, 13288
Marseille-Luminy, and \\ PHYMAT, Universit\' e de Toulon et du Var,
83957 La Garde, France}
\authors{P. Exner}
\address{Nuclear Physics Institute, AS, 25068 \v Re\v z near Prague,
and \\ Doppler Institute, Czech Technical University, 11519 Prague,
Czech Republic}
\authors{and B. Meller\footnote{Address as of Oct 1st 1997: 
Facultad de F{\'\i}sica, P.U. Cat{\'o}lica de Chile, 
Casilla 306, Santiago 22, Chile}}
\address{Centre de Physique Th\'eorique, CNRS, 13288
Marseille-Luminy, and \\ PHYMAT, Universit\' e de Toulon et du Var,
83957 La Garde, France}
\abstract{We consider curvature--induced resonances in a planar
two--dimensional Dirichlet tube of a width $\,d\,$. It is shown that
the distances of the corresponding resonance poles from the real axis
are exponentially small as $\,d\to 0+\,$, provided the curvature of
the strip axis satisfies certain analyticity and decay requirements.}

\section{Introduction}

Spectral and scattering properties of Dirichlet Laplacians in curved
tubes have attracted a wave of physical interest attention recently,
because they provide models of some quantum systems which new
experimental techniques made possible to construct, such as
semiconductor quantum wires --- see, \eg,
\cite{ABGM,Ba,CLMM,DE,Sa,SRW,VOK1} and references therein --- or
hollow--fiber atomic waveguides \cite{SMZ}, and because they exhibit
some unexpected mathematical properties leading to new physical
effects.

The key observation is that a nonzero curvature gives rise to an
effective interaction which produces localized solutions of the
corresponding Helmholtz (or stationary Schr\"odinger) equation ---
\cf\cite{ES,GJ} and the review paper \cite{DE} --- with eigenvalues
below the bottom of the continuous spectrum. The same mechanism is
responsible for a nontrivial structure of the scattering matrix
manifested by resonances in the vicinity of all the higher
thresholds. These resonances modify substantially transport
properties of such a ``quantum waveguide''; they have been observed
in numerically solved examples, for instance, in
\cite{SM,VKO,VOK1,VOK2,WS}.

On the mathematical side, it was shown in \cite{DES} that if
a curved planar strip has a constant width $\,d\,$ which is small
enough, and if the strip--axis curvature satisfies certain regularity
and analyticity assumptions, there is a finite number of resonances
in the vicinity of the higher thresholds (which coincides with the
number of isolated eigenvalues below the bottom of the continuous
spectrum). Moreover, an expansion of the resonance--pole positions in
terms of $\,d\,$ was derived and the imaginary part of the first
non--real term given by the ``Fermi golden rule'' was shown to be
exponentially small as $\,d\to 0+\,$.

The present paper addresses the question whether also the 
{\em total resonance widths\/}
are exponentially small as $\,d\to 0+\,$.
We give an affirmative answer under
essentially the same assumptions as used in \cite{DES} 
and obtain the same expression for the exponential factor in the bound.
Furthermore the
exponential factor we obtain coincides with the heuristic
semi-classical prediction, \cf \cite{LL,DES}.

As explained in \cite{DES} and in section~2.3 below, the mechanism
behind the formation of these resonances is a tunneling effect, however, in
the ``momentum direction''. To estimate this effect we therefore need
exponential bounds on eigensolutions in the Fourier representation.
This is a novelty and a difficulty, since the Schr\"odinger operator
becomes nonlocal in this representation. To deal with such nonlocal
operators we have developed an appropriate functional calculus based
on the Dunford--Taylor integral. This made possible, in particular, to
extend to such a situation the complex deformations of operators and
the Agmon method \cite{Ag}. These new techniques has already been
announced in \cite{DM1}. Here they are presented in detail for the
case of bounded nonlocal operators; an extension to some unbounded
cases will be given in \cite{DM2}.

Rigorous analysis of tunneling in phase space is a rather new field of
interest. Some recent works on this topic based on micro-local
analysis and pseudo-differential techniques can be found in
\cite{HeSj2,Ma,N}. In particular L.Nedelec \cite{Ne} has recently
obtained our Theorem~2.2 with such methods.
\setcounter{equation}{0}
\section{The results}

\subsection{Preliminaries}

Let us recall briefly the problem; for more details we refer to
\cite{DES}. The object of our interest is the Dirichlet Laplacian
$\,-\Delta_D^{\Omega}\,$ for a curved strip $\Omega\subset \R^2$ of a
fixed width $\,d\,$. We exclude the trivial case and make a global
restriction:
      \begin{description}
   \item{(a0)} $\;\Omega\,$ is not straight and does not intersect
itself.
      \end{description}
If the boundary of $\,\Omega\,$ is sufficiently smooth --- which
will be the case with the assumptions mentioned below ---  one can
check using natural curvilinear coordinates that
$\,-\Delta_D^{\Omega}\,$ is unitarily equivalent to the
Schr\"odinger type operator
      \begin{equation} \label{Hamiltonian}
H\,:=\,-\partial_s b\partial_s -\partial_u^2 +V
      \end{equation}
on the Hilbert state space on the ``straightened'' strip, $\,\HH:=
L^2(\R\!\times\!(0,d), ds\,du)\,$, with the Dirichlet condition at
the boundary, $\,u=0,d\,$, where $\,b,\,V\,$ are operators of
multiplication by the functions
      \begin{eqnarray} \label{kinetic factor}
b &\!:=\!& (1\!+\!u\gamma)^{-2}\,, \\ \nonumber \\
\label{potential}
V &\!:=\!& -\,{1\over 4}\, b\gamma^2 +\,{1\over 2}\, b^{3/2} u\gamma''
-\,{5\over 4}\, b^2 u^2\gamma'^2\,,
      \end{eqnarray}
respectively, and the function $\,\gamma:\,\R\to\R\,$ in these
relations characterizes the strip boundary $\,u=0\,$ through its
signed curvature $\,\gamma(s)\,$ at the point tagged by the
longitudinal coordinate $\,s\,$ --- \cf\cite{ES,DE}.

Let us list now the used hypotheses. In addition to the assumption
(a0), we shall suppose that
      \begin{description}
   \item{(a1)} $\;\gamma\,$ extends to an analytic function in
$\,\Sigma_{\alpha_0,\eta_0}:= \{\, z\in\C\,:\; |\arg(\pm z)|<\alpha_0
 \;$ or
$\,|\im z|<\eta_0\,\}\,$ with $\,\alpha_0<{\pi\over 2}$
and $\,0<\eta_0\,$;
for the sake of simplicity we denote it by
the same symbol.
   \item{(a2)} For all $\,\alpha<\alpha_0\,$ and
all $\,\eta<\eta_0\,$ there are positive constants $\,c_{\alpha,\eta}\,$
and $\,\eps\,$
such that
$\;|\gamma(z)|<c_{\alpha,\eta}(1\!+\!|z|)^{-1-\eps}\,$ holds in
$\,\Sigma_{\alpha,\eta}\,$.
       \end{description}
By an easy application of the Cauchy formula, the assumptions (a1)
and (a2) imply that the derivatives of $\,\gamma\,$ satisfy
$$
|\gamma^{(r)}(z)|\,<\, c_{r,\alpha,\eta} (1\!+\!|z|)^{-1-r-\eps}
$$
in $\,\Sigma_{\alpha,\eta}\,$ for any $\,\alpha<\alpha_0\,$ and
any $\,\eta<\eta_0\,$.
This yields for the potential (\ref{potential})
the bound
      \begin{equation} \label{potential bound}
|V(z,u)|\,<\, c'_{\alpha,\eta} (1\!+\!|z|)^{-2-\eps}
      \end{equation}
with some $\,c'_{\alpha,\eta}>0\,$ for all $d$ small enough.

We are interested in resonances of $\,H\,$ which are understood in
the standard way \cite{AC,RS,Hu}: suppose that the function
$\,z\mapsto F_\psi(z):=((H-z)^{-1}\psi,\psi)\,$ admits a meromorphic
continuation from the open upper complex half-plane to a domain in the
lower half-plane for $\,\psi\,$ from a dense subset
$\,\AAA\subset\HH\,$. If $\,F_{\psi}$ has a pole for some
$\,\psi\in\AAA\,$, we call the former a {\em resonance\/}.

Resonances are often obtained as perturbations of an operator with
eigenvalues embedded in the continuous spectrum. This is also the
case in our present situation; the corresponding comparison operator
is
      \begin{equation} \label{comparison operator}
H^0\,:=\,A-\partial^2_u\,, \qquad A\,:=\,-\partial^2_s+V^0\,,
      \end{equation}
with $\,V^0(s):=V(s,0)=\! -{1\over 4}\, \gamma(s)^2\,$ and domain
$\DD(H^0):=\HH^2(\R)\otimes (\HH_0^{1}\cap
\HH^{2})((0,d))$, where $\,\HH^{n}$
and $\,\HH^n_0$ are the usual Sobolev spaces.
The perturbation is defined by
      \begin{equation} \label{perturbation}
W\,:=\,H-H^0\,.
      \end{equation}
The spectrum of the operator $\,H^0\,$ is of the form
$$
\sigma(H^0)\,=\,\left\lbrace\, \lambda\!+\!E\,:\; \lambda\in\sigma(A),\;
E\in\sigma(-\partial_u^2)\,\right\rbrace\,,
$$
where
$$
\sigma(A)\,=\, \left\lbrace\, \lambda_n\, \right\rbrace_{n=1}^N\,\cup
\,[0,\infty)\,, \qquad
\sigma(-\partial_u^2)\,=\, \left\lbrace\, E_j\,
\right\rbrace_{n=1}^{\infty}
$$
with $\,E_j:= \left({\pi j\over d}\right)^2$. Since $\,\int_{\R}
V^0(s)\, ds<0\,$ due to (a0),
the discrete spectrum of $\,A\,$ is nonempty. The
eigenvalues $\,\lambda_n\,$ are negative, simple, and their number
$\,N\,$ is finite in view of the bound (\ref{potential bound}) ---
\cf\cite[Sec.XIII.3]{RS}. Then the eigenvalues
$$
E^0_{j,n}\,:=\, \lambda_n+E_j
$$
above $\,E_1\,$ are embedded in the continuous spectrum of $\,H^0\;$;
for small enough $\,d\,$ this occurs for all $\,j\ge 2\,$ and
$\,n=1,\dots,N\,$.

An alternative way to express the operator $\,H^0\,$ and functions of
it, is through the transverse--mode decomposition. Denote by
      \begin{equation} \label{transmode}
\chi_j(u)\,:=\,\sqrt{2\over d}\,\sin\left({\pi ju\over d}\right)\,,
\qquad j=1,2,\dots,
   \end{equation}
the eigenfunctions of $\,-\partial_u^2\,$ corresponding to the
eigenvalues $\,E_j\,$. Let $\,\JJ_j\,$ be the embedding
$\,L^2(\R,ds)\,\to\,L^2(\R,ds) \otimes\chi_j\,\subset\,\HH\;$; the
adjoint of this operator is $\,\JJ_j^*:\, \HH\to L^2(\R,ds)\,$ acting
as $\,(\JJ^*_jf)(s)= \int_0^d f(s,u)\chi_j(u)\,du\,$. Given
$\,j\in\N\,$, we denote by $\,P_{j}\,$ the projection onto
the mode with index $\,j\,$, $\;P_{j}:= \JJ_j\JJ_j^*\,$,
 and set $\,Q_{j}:= I_{\HH}\!-P_{j}\,$.

The perturbation $\,W\,$ consists of operators which
 couple  different transverse modes. As a
result the embedded eigenvalues turn into resonances. With our
assumptions the result of \cite{DES} holds:

      \begin{theorem} \label{resonances}
Assume (a0)--(a2). For all sufficiently small $\,d\,$ each eigenvalue
$\,E^0_{j,n}\,$ of $\,H^0\,,\; j\ge 2\,,n=1,\dots,N\,$, gives rise to
a resonance 
$\,E_{j,n}(d)\,$ of $\,H$, the position of which is given by a
convergent series 
      \begin{equation} \label{resonance expansion}
E_{j,n}(d)\,=\, E_{j,n}^0+ \sum_{m=1}^{\infty} e_m^{(j,n)}(d)\,,
      \end{equation}
where $\,e_m^{(j,n)}(d)= \OO(d^m)\,$ as $\,d\to 0+\,$. The first term
of the series is real--valued, and the second satisfies the bound
      \begin{equation} \label{Fermi rule bound}
0\,\le\, -\im e_2^{(j,n)}(d)\,\le\, c_{\eta,j}\, e^{-2\pi\eta
 \sqrt{2j-1}/d}
      \end{equation}
for any $\,\eta\in(0,\eta_0)\,$, the constant $\,c_{\eta,j}\,$
depending on $\,\eta\,$ and $\,j\,$.
      \end{theorem}
\subsection{Main theorems}
Our aim in this paper is to show that similar bounds can be proven
for the total resonance width. This is the contents of the following
two theorems:

      \begin{theorem} \label{full width bound}
Assume (a0)--(a2). Then for any $\,\eta\in(0,\eta_0)\,$, $j\,\geq
2\,$ and $n=1,\dots,N\,$ there is
$\,C_{\eta,j}>0\,$ such that
      \begin{equation} \label{resonance width bound}
0\,\le\, -\im E_{j,n}(d)\,\le\, C_{\eta,j}\, e^{-2\pi\eta
 \sqrt{2j-1}/d}
      \end{equation}
holds for all $\,d\,$ small enough.
      \end{theorem}

      \begin{theorem} \label{instanton resonance bound}
In addition, assume that $\,\gamma\,$ extends to a meromorphic
function in $\,\Sigma_{\alpha_0,\eta_1}\,$ with $\,\eta_1>\eta_0\,$.
Let $\,\eta_p<\eta_1\,$ be the minimal distance to the real axis of the
poles
and assume that the maximal order of the poles at this distance is
$\,1\leq m < \infty\,$; then there are positive constants $\,C^{(1)}_{j}\,$ and
$\;C^{(2)}_{j}\,$ such that
      \begin{equation} \label{instanton bound}
0\,\leq\,-\im E_{j,n}(d)
\le C^{(1)}_{j}\,
\exp\left\lbrace\, -{2\pi\eta_p\over d}\, \sqrt{2j\!-\!1}\,
\left(1\!-\!C^{(2)}_{j}\,d^{1/(m+1)}\,\right)
\right\rbrace \phantom{AAA}
      \end{equation}
holds for all $\,d\,$ small enough.
      \end{theorem}
 \begin{remarks}
{\em
(i) There is an heuristic prediction for the value of $\im E_{j,n}(d)$ based 
on a formal semi-classical analysis where the role of the semi-classical 
parameter is played by $d$ as $d$ tends to zero. What one expects according 
to this prediction (for the details we refer the reader to \cite{DES}, in 
particular Remark~4.2e therein and also to the scheme of the proof below) is 
that $\im E_{j,n}(d)$ should behave like 
$\displaystyle C_j(d)\exp\left(-\frac{2\pi\eta_0}{d} \sqrt{2j-1}\right)$ 
where $C_j(d)$ is 
polynomially bounded in $d^{-1}$.  However there is no chance to get such a 
precise behaviour without knowing the type of singularity that the curvature 
$\gamma$ exhibits at distance $\eta_0$ from the real axis.
This is why in Theorem~2.2 we lose an arbitrary small 
part of the exponential decay rate 
and get a prefactor $C_{j,\eta}$ which may eventually 
diverge as $\,\eta\,$ tends to $\,\eta_0\,$. 
This kind of bound is typical in such 
a semi-classical context, see e.g \cite{Ag}.

(ii) The merit of Theorem~2.3 is to show that with a precise assumption
on the type of singularity of $\gamma$ we are able to produce a bound which 
has a leading behaviour in accordance with the  heuristic prediction. We would
like to stress that, to our knowledge, this is the most precise bound obtained
so far on the total resonance width in such a situation.
}
\end{remarks}

Since we shall deal in the following mostly with a single resonance,
we drop the subscripts $\,j,n\,$ as well as the argument $\,{d}\,$ 
whenever they
are clear from the context.
\subsection{A sketch of the proofs}
Consider first the system described by the decoupled Hamiltonian
$H^0$. Each state $\,\phi\,$ of this system can be decomposed into the sum
of its transverse modes, $\phi^j\otimes\chi_j, j=1,2,\dots\,$, and
this decomposition is invariant under the dynamics generated by
$H^0$.  For each channel $j$ the dynamics of $\phi^j$ is governed by
the ``longitudinal'' Hamiltonian $\, A+E_j =-\partial ^2_s + V^0
+E_j\,$ which, due to the nonzero curvature of the guide and its
decay at infinity (by (a0) and (a2)), possesses either bound states
for energies below  $\,E_j\,$ or scattering states otherwise. Fix now
a $\,j\ge 2\,$ and suppose that $\,d\,$ is small enough so that a
given bound--state energy $\,E^0= E_j\!+\!\lambda_n\,$ of
$\,A\!+\!E_j\,$ is embedded in the continuous spectrum of the lower
modes. The only possible solutions to the equation
$\,H^0\phi=E^0\phi\,$ are then this bound state in the $\,j$--th mode
and $\,j\!-\!1\,$ scattering states in the modes below. This
structure is reflected in the classical phase space portrait of $H^0$
at energy $E^0$ (see Figure~1; we recall that for a matrix
Hamiltonian $\,H^{(cl)}(q,p)\,$ the energy shell at $\,E^0$ is given by
$\,\det (H^{(cl)}(q,p)-E^0)=0\,$); the energy shell of $H^0$ is the
union of the curves $\,p^{(cl)}_k(s):=\pm\sqrt{E^0\!- E_k- V^0(s)},\
k=1,\dots ,j\,$. As expected only the $\,j$--th curve is compact.
\begin{figure}
\epsfig{file=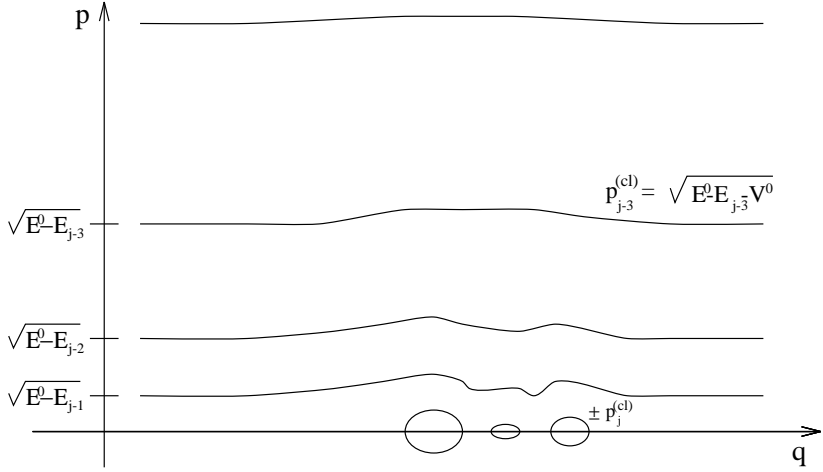} 
\caption{A schematic phase-space portrait of a bent waveguide}  
\end{figure}

Let now $\phi$ be a bound state of $\,H^0$ in the $\,j$--th
channel and consider its evolution under the full dynamics given by
$\,H=H^0+W$. In general, various channels of $\,H^0$ are now coupled
by $W$ and $\phi$ will undergo transitions to all other energetically
allowed channels. For $d$ small enough there will be no significant
changes of the classical phase space portrait by the addition of $W$.
Thus for the classical dynamics no transition is possible between 
different channels. Hence the transitions in the quantum dynamics are
of the tunneling type, but in contrast to the usual situation the
tunneling takes place in the $\,p\,$ (\ie momentum) direction.
More precisely, the projection of the energy shell on the $p$--axis
consists of intervals of classically allowed momenta, one for the
$\,j$--th mode situated at the origin and two for all other modes with
index $k<j$, situated around $\,\pm\sqrt{E^0-E^k}\,$. They are
separated by ``gaps'', \ie, classically forbidden (momentum) regions
which have a size of order $d^{-1}$ as $\,d\,$ tends to zero. The
existence of such gaps suggests  that the solutions of $H\phi=E\phi$,
with  $E$  close to $E^0$,
decay exponentially as functions of $\,p\,$ in these gaps, a key property in
the sequel. The main contribution to this tunneling process should
come from the transitions from $p^{(cl)}_j$ to $p^{(cl)}_{j-1}$, since it
is the first gap that the quantum state has to cross. This motivates
our decomposition of the momentum space into
$\overline{\Omega}_i\cup\Omega_e$, with
$\Omega_i\cap\Omega_e=\emptyset$ and $\Omega_i:=(-\omega,\omega)$,
where $\omega$ is approximately equal to $\sqrt{E_j-E_{j-1}}$, the
size of the first gap.

Let $\,H_{\theta}\phi_{\theta}=E\phi_{\theta}\,$ be the eigenvalue
equation for a resonance $\,E\,$ and the corresponding resonance function
$\,\phi_{\theta}$ and assume that this resonance is associated to
$E^0_{j,n}=E^0$ as in Theorem~2.1. The complex deformation, denoted by
$\theta$, is chosen as a scaling of the momentum exterior to 
$\,\Omega_i$, \cf (3.1).
In the transverse mode decomposition the
eigenvalue equation becomes an infinite system of coupled
differential equations in $L^2(\R)$ which can be solved for the
$\,j$--th component $\phi_{\theta}^j:= \JJ^*_j\phi_{\theta}$ of
$\phi_{\theta}$ leading to the following equation
\begin{equation}
\label{2.12}
E\phi_{\theta}^j \,=\,(H^j_\theta-B^j_\theta(E))\phi_{\theta}^j
\end{equation}
in $\,L^2(\R)\,$, where $\,H^j_{\theta}:= \JJ^*_j
H_{\theta}\JJ_j,\;B^j_{\theta}(E):=
\JJ^*_j W_\theta \widehat{R}^j_{\theta}(E) W_\theta\JJ_j\;$ and
$\,\widehat{R}^j_{\theta}(E):=Q_j(Q_j H_{\theta}Q_j\!-E)^{-1}Q_j\,$.
In section~4, with the help of (\ref{2.12}) we prove stability of the
spectral value $E^0$ of $H^0_\theta$ under the perturbation by
$W_\theta$. Then we are able to show that for $d$ small enough the
tunneling picture given in the previous paragraph is correct. Indeed
we obtain the following exponential bound on the $\,j$--th component of
$\phi_\theta$: let $\rho$ be a function obeying (3.2), then
\begin{equation}
\label{2.13}
\left\|-i\partial_s\,e^{\rho(-i\partial_s)}\phi_{\theta}^j \right\|^2 +
\left\|e^{\rho(-i\partial_s)}\phi_{\theta}^j \right\|^2 < \infty\,,
\end{equation}
the bound being uniform as $d$ tends to zero. This is one of the main
ingredients of this paper. We turn now to the explanation of how one can
use (\ref{2.13}) to derive our exponential estimate on $\im E$, which
is the purpose of Section~6.

From (\ref{2.12}) we obtain by straightforward algebraic computations
the following equation for $\im E$:
\begin{eqnarray}\label{2.14}
\im E &\!=\!&
((\im H^j_\theta-Z_\theta(E))\phi_{\theta}^j,\phi_{\theta}^j)\\
\nonumber
Z_\theta(E)
&\!:=\!&
\JJ^*_j\left\{2\re [\im (W_\theta)\widehat{R}^j_\theta W_\theta]-
(\widehat{R}^j_\theta
W_\theta)^*\,\im\widehat{H}^j_\theta\,\widehat{R}^j_\theta W_\theta
\right\}\JJ_j
\end{eqnarray}
with $\widehat{H}^j_\theta:=Q_jH_\theta Q_j\,$ provided $\phi_\theta$
is chosen with the unit norm. The merit of this cumbersome formula is
that it indicates that the operator $\im H^j_\theta-Z_\theta(E)$
should acts as a localization on $\Omega_e$ in the momentum space. This
is due to the fact that each of its three terms contains an imaginary
part of a scaled operator which is expected to vanish on $\Omega_i$
where the deformation does not operate. If this localization property
would be true then (\ref{2.14}) combined with (\ref{2.13}) would give
immediately the desired exponential estimate on $\im E$: 
$$ 
|\im E|\leq {\rm const\:} e^{-2\rho(\omega)}\,.  
$$ 
Unfortunately, since most of the operators involved here are
non-local in the momentum variable, this simple reasoning does not work.
However we are able to show directly that this localization property
is valid in the following weaker sense 
$$ 
e^{-\rho}|\im
H^j_\theta-Z_\theta(E)|e^{-\rho}\leq {\rm const\:}
e^{-2\rho(\omega)}(-\partial_s^2 +1) 
$$ 
which is all what we need.

Let us finish the survey of the paper contents.
To deal with the Schr\"odinger operator in the momentum representation,
and in particular, with its image under an exterior scaling in the
momentum variable, we have developed in Section~3
a functional calculus  based on the Dunford-Taylor integral.
All the necessary material for
the exterior scaling is presented in Section~3.
Finally the
extension of our method to the case where the nearest singularity of
the curvature in the complex plane is a pole is done in Section~7.
\setcounter{equation}{0}
\section{Complex scaling and functional calculus}                  
From this moment on we pass to the unitary equivalent situation by
performing the  inverse Fourier transformation in the $\,s\,$ variable,
denoted by $\,F_s^{-1}\,$. We introduce the notation:
$$
p\,:=\,F_s^{-1}\,i\partial_s \,F_s \qquad \mbox{and}\qquad 
D\,:=\, -i\partial_p\,=\,F_s^{-1}\,s \,F_s 
$$
For
all other transformed operators we shall use 
the same symbols as before:
$$
H=p\,b(D,u)\,p -\partial_u^2 +V(D,u).
$$
Note that now $\,\DD(H^0):=\DD(p^2)\otimes (\HH_0^{1}\cap
\HH^{2})((0,d))\,$.
\subsection{Exterior scaling in momentum representation}
Complex dilations represent a useful tool to reveal resonances in
systems with Hamiltonians having certain analytic properties. In the
present case, we use the exterior dilation defined as follows:
      \begin{equation} \label{dilation}
p_{\theta}(t)\,:=\, \left\lbrace\;    \begin{array}{lll} t & \quad \mbox{if}
\quad & t\in\Omega_i:= (-\omega,\omega) \\ \\ \pm\omega+
e^{\theta}(t\!\mp\!\omega) & \quad \mbox{if} \quad & 
t\in\Omega_e:=\R\setminus \overline{\Omega}_i
   \end{array} \right.
      \end{equation}
where $\,\omega\,$ is a  positive number to be determined later.
The parameter
$\,\theta\,$ 
 takes complex values in a strip around the real axis;
defining the sets $\,\SS_\alpha :=\{\theta\in\C,\,
|\im\theta|<\alpha\}\,,\alpha>0$ we have
$\,\theta\in\SS_{\alpha_0} $. The function
$\,p_\theta\, $ is for real $\,\theta\,$  a
piecewise differentiable homeomorphism of $\R$ whose
unitary implementation on $\,L^2(\R)\,$ is defined by
$$
U_\theta\varphi := {p'_\theta} ^{1/2}\varphi\circ p_\theta.
$$
$U_\theta\,$ and $\,p_\theta\,$ are both called (exterior) dilation.
In general, to denote the image under this dilation we use the
index $\,\theta\,$. 

Recall how one uses $\,U_\theta\,$
to get a complex deformation of operators. With a given closed operator
$\,T\,$, one constructs the  family of operators for
$\,\theta\in\R\,$:
$$
\theta\to T_\theta := U_\theta T\,U_\theta^{-1}.
$$
If this function has an analytic extension to some $\,\SS_\alpha $ (in a suitable
sense
--- \cf \cite[Ch.VII]{Ka}), the resulting family is what is usually
called a complex (family of) deformation(s) of $\,T\,$.

We begin by considering the complex deformation of $p$ and $D$:

      \begin{proposition} \label{p_theta}
(i) $\;\{\,p_{\theta}^2:\; \theta\in\C\,\}\,$ is a self--adjoint
family of type A in the sense of \cite{Ka}
with common domain $\,\DD(p_{\theta}^2)= \DD(p^2)\,$. \\
(ii) $\;\sigma(p_{\theta}^2)= [0,\omega^2] \cup p_\theta^2(\Omega_e)\,$.
      \end{proposition}
We would like to remark here that since we are scaling in the Fourier
image, we will have to use $\,\theta\/$'s
with a negative imaginary part to
make the essential spectrum turn into the lower complex half-plane.
      \begin{proposition} \label{D_theta}
(i) $\;\{\,D_{\theta}\,:\; \theta\in\C\,\}\,$ is a self--adjoint
analytic family. A vector $\,u\,$ belongs to $\,\DD(D_{\theta})\,$ \Iff
$\,u\in \HH^{1}(\Omega_i)\oplus \HH^{1}(\Omega_e)\,$ and
$\,u(\pm\omega\pm 0)= e^{\theta/2}u(\pm\omega\mp 0)\,$, the action
of the operator being given by
$$
(D_{\theta}u)(t)\,=\, (p'_{\theta})^{-1}(D\,u)(t)\,=\, \left\lbrace
   \begin{array}{lll} -iu'(t) & \quad \mbox{if}
\quad & t\in\Omega_i \\ \\ -ie^{-\theta}u'(t) & \quad \mbox{if} \quad &
t\in\Omega_e    \end{array} \right.
$$
(ii) $\;\sigma(D_\theta)= e^{-\theta}\, \R\,$.
      \end{proposition}
{\em Proof:} (i) By the standard argument --- \cf \cite{CDKS} for the
case of the Laplacian. \\
(ii) This is a straightforward calculation using the
explicit expression of the resolvent kernel of
$D_\theta$ and bounding it by the Schur-Holmgren norm.
Recall that the Schur-Holmgren norm of an integral operator $\,B\,$
is defined through its kernel as
$$
\|B\|_{SH}\,:=\, \max\left\lbrace\, \sup_x \int |B(x,y)|\,dy\,,\;
\sup_y \int |B(x,y)|\,dx \right\rbrace\,.
$$
(\cf\cite[Example III.2.4]{Ka}) and majorizes the operator norm.
\quad\QED\\
\vspace{3pt}\noindent
To study a complex deformation of operators of the form
$f(D)$ we need to develop a functional calculus for $\,D_{\theta}\,$.
\subsection{Functional calculus for $\,D_{\theta}\,$}
Since the operator under consideration contains
the metric term (\ref{kinetic factor}) and the potential
(\ref{potential}), we have to define the corresponding
operator functions. The standard functional calculus is not
applicable here, because $\,D_{\theta}\,$ is not even normal
for complex $\,\theta\,$; instead we use the Dunford-Taylor
integral.
The original theory for unbounded operators is
 exposed in \cite[Sec.VII.9]{DS}.
But since analytic families of operators are not treated there and
since it is necessary for our estimates to modify the original
definition, we present the adapted theory in detail.

   \begin{definition} \label{DTI}
Let a function $\,f:\C\to\C\,$ satisfy the same requirements as
$\,\gamma\,$ in (a1) and (a2). Suppose that $\,T\,$ is a closed
operator in $\,L^2(\R)\,$ and there is an open set $\,\VV\,$ which
obeys strict 
inclusions $\,\sigma(T)\subset \VV\subset \Sigma_{\alpha_0,
\eta_0}\,$ and whose boundary $\,\partial \VV\,$ consists
of a finite number of rectifiable Jordan curves with a positive
orientation. Suppose also that
$\,(T\!-\!z)^{-1}$ is uniformly bounded on $\,\partial \VV\,$.
Then we define
$$
f(T)\,:=\, {i\over 2\pi}\, \int_{\partial \VV}\, f(z) (T\!-\!z)^{-1}\, dz\,.
$$
The operators defined this way will be called Dunford-Taylor operators.
   \end{definition}

   \begin{lemma} \label{DTI properties}
(i) $\;f(T)\,$ is a well--defined bounded operator on $\,L^2(\R)\,$. \\
(ii) If $\,T\,$ is self--adjoint, $\,f(T)\,$ coincides with the
operator obtained by the usual functional calculus. \\
(iii) $\;Bf(T)B^{-1}= f(BTB^{-1})\,$ holds for any bounded
operator $\,B\,$ with bounded inverse. \\
(iv) For $\,\theta\in\R\,$, let $\,T_{\theta}:= U_{\theta}T
U_{\theta}^{-1}\,$
such that
$\,\{T_{\theta},\theta\in\SS_\alpha\}\,,0<\alpha\leq\alpha_0\,$
 is an analytic family of
operators.
Assume that there is an open set $\,\VV\,$ with
$\displaystyle\,\bigcup_{\theta\in\SS_\alpha }\!
\sigma(T_{\theta}) \subset \VV \subset \Sigma_{\alpha_0,
\eta_0}\,$, that $\,(T_{\theta}-z)^{-1}\,$ is uniformly bounded on
$\,\partial \VV\,$ for all $\,\theta\in \SS_\alpha$
and that $\,\partial \VV\,$ obeys the conditions of the definition. Then
for all $\,\theta\in \SS_\alpha$
$$
\Big(f(T)\Big)_{\theta}\,=\, f(T_{\theta})
$$
and these operators form a bounded  analytic
family on $\,L^2(\R)\,$.
  \end{lemma}
{\em Proof:\/} We prove this lemma in appendix A.
\vskip3pt
In our case $\,T:=D_{\theta}\;$;
since $\,\sigma(D_{\theta})=e^{-\theta}\R\,$, it is clear that for any
$\,\theta\,$ with $\,|\im\theta|<\alpha_0\,$ there is a
domain
$\,\VV_\theta\,$
satisfying strict inclusions
$\,\sigma(D_{\theta})\subset \VV_\theta
\subset \Sigma_{\alpha_0,0}\,$.
But we still need to control the resolvent of $\, D_\theta\,$.

Before doing that we want to introduce another operator deformation
we shall need later:
$$
T_\rho \, := \, e^\rho\,T\,e^{-\rho},
$$
is usually called the image of $\,T\,$ under
the {\em boost\/} $\,-i\rho\,$, where $\,\rho\,$ is an
absolutely continuous function.
In particular, $\, D_\rho = D + i\rho'$, suggesting the origin of this
terminology. We shall only consider real functions for the boosting,
\ie\hskip3pt only purely imaginary boosts. Furthermore it will be sufficient
 for our purpose to use only boosts being constant on $\, \Omega_e\,$.
Then, of course, the boost commutes  with the exterior scaling and
there should be no confusion concerning  our notation,  $\,
T_{\theta,\rho}\,$, for the indication of the two deformations of $\,T$.

If $\,\rho'\,$ is supported on $\,\Omega_i\,$, one has
$\,D_{\theta,\rho}=D_{\theta} + i\rho'\,$.
Note that we write  $\,D^*_{\theta,\rho}\,$ for
$\,(D_{\theta,\rho})^*\,$. 
Finally since $\rho$, and therefore also $\,e^{\pm\rho}$ are bounded,
we have  $\,e^\rho f(T)e^{-\rho}=f(T_{\rho})\,$
by Lemma~\ref{DTI properties} (iii).

Due to group property of the exterior dilation in $\,\theta\,$ it is
sufficient to perform the estimates for purely imaginary
$\,\theta\,$; we shall write and prove the corresponding bounds for
$\,\theta=i\beta\,,\;\beta\in\R\,$ only.

We will use the symbol
$\,\chi_{A}\,$ to denote the
characteristic function of a  set $\,A\,$.
      \begin{proposition} \label{general reso p estimate}
Let $\rho$ be a real, bounded, absolutely continuous function
on $\R$ which is
constant on $\Omega_e$. Furthermore define
$\,\GG_{\beta,\rho} :=\{z\in\C :
|\arg\left(\pm ze^{-i\beta/2}\right)|\leq {|\beta|/2}\;
{\rm or}\;
|\im z|\leq \|\rho'\|_\infty \}
\,$ --- \cf Figure~2.
Then for all $\beta$ with $\,|\beta|<{\pi\over 2}\,$ and all
$\,z\in\C\setminus\GG_{\beta,\rho}\,$
$$
\left\|\,(D_{i\beta,\rho}\!-\!z)^{-1}\right\|
\leq \,{\rm dist}(z,\GG_{\beta,\rho})^{-1}
$$
     \end{proposition}
{\em Proof:} Let $\,v\in \DD(D_{i\beta})\,$ and $\,w:=p'_{-i\beta}v\,$.
Then $\,\|w\|\leq\|v\|\,$;
we have
   \begin{eqnarray*}
\left\|\,(D_{i\beta}\!+\!i\rho'\!-\!z)v\right\|\|w\|
 &  \!\geq\!&
|(p'_{i\beta}(z\!-\!i\rho'\!-\!D_{i\beta} )v,v)|
\,\geq\, |\im(p'_{i\beta}(z \!-\!i\rho'\!-\!D_{i\beta} )v,v)| \\ \\
 &\!=\!&
|((\im p'_{i\beta}z \!-\!\rho')v,v)|\,;
   \end{eqnarray*}
The last equality is due to:
   \begin{eqnarray*}
\im(p'_{i\beta}\,D_{i\beta}v,v)
&\!=\!&{1\over 2}\left\{(\chi_{\Omega_i}v',v)+(\chi_{\Omega_i}v,v')+
(\chi_{\Omega_e}v',v)+(\chi_{\Omega_e}v,v')\right\}\\ \\
&\!=\!&
{1\over 2}\left\{ |v|^2\right|^{+\omega-0}_{-\omega+0}\left.\left.
+|v|^2\right|^{-\omega-0}_{+\omega+0}\right\}
\\ \\
&\!=\!& 0\,,
    \end{eqnarray*}
using that for $\,v\in \DD(D_{i\beta})\,$ the discontinuity at
$\,\pm\omega\,$
is just the phase $\, e^{i\beta/2}$.
\quad  \QED
\vskip3pt

For the sake
of brevity we shall use the shorthand $\,f_\vartheta:=f(D_\vartheta)\,$
 for the
Dunford--Taylor operators under consideration, where $\,\vartheta=
\theta,\,\rho\,$, {\em etc.;} the superscript $\,c\,$ will denote the
complement of a set. The last proposition implies, in particular, that
$$
\max\left\{\,\left\|\,(D_{i\beta,\rho}\!-\!z)^{-1}\right\|,\,
\left\|\,(D^*_{i\beta,\rho} \!-\!z)^{-1}\right\|\, \right\}
\,\leq \,\mbox{dist}(z,\Sigma_{|\beta|,\|\rho'\|_\infty})^{-1}
$$
holds for $\,z\in (\overline\Sigma_{|\beta|,\|\rho'\|_\infty})^{c}\,$; the
involved
domains are sketched on Figure~2. Furthermore we want to state the
general conditions on $\,\rho\,$ which will be imposed up to the end of
Section~6:
\newlength{\mpwidth}
\setlength{\mpwidth}{\textwidth}
\addtolength{\mpwidth}{-30truemm}
\begin{equation}\label{condition on rho}
\left\{
        \begin{array}{ll}\mbox{(i)}&
\mbox{\parbox[t]{\mpwidth}{
 $\,\rho\,$ is a real, absolutely continuous function on $\R$ which
is constant on $\,\Omega_e\,$, }}\\
\mbox{(ii)} &\|\rho'\|_\infty \leq \eta <\eta_0.
        \end{array}\right.
\end{equation}
The functions $b$ and $V$ --- \cf (\ref{kinetic
factor}) and (\ref{potential}) --- are understood as rational functions of
$u\gamma$,
$u\gamma'$ and
$u\gamma''$. Notice that their structure is particularly simple; there are
only powers of $1+u\gamma$ appearing in the denominator. 
\begin{figure}
\epsfig{file=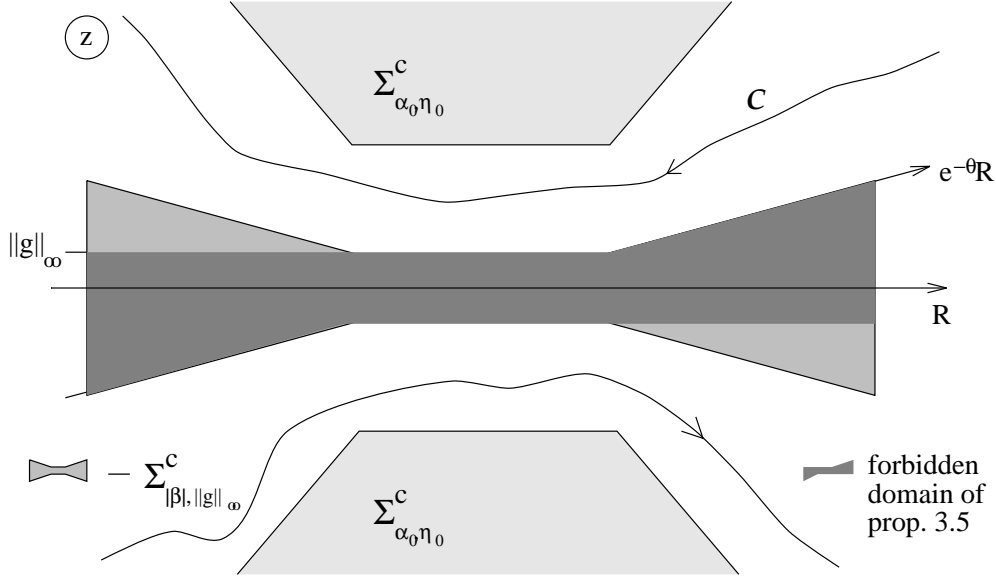} %
\caption{To the definition of the Dunford-Taylor integral
(the case $\im\theta \!<\! 0$)}
  \end{figure}
        \begin{proposition} \label{DTI gamma}
Let $\,\rho\,$ 
satisfy (\ref{condition on rho}).\\
(i)  
$\,\gamma_{\theta,\rho}\,$ and therefore $\,V^0_{\theta,\rho}\,$, as
well as $\,\gamma'_{\theta,\rho}\,$ and $\,\gamma''_{\theta,\rho}\,$
are bounded  self--adjoint analytic families
of operators in $L^2(\R)$ in $\theta$ provided
$\,\theta\in\SS_{\alpha_0}\,$ .\\  
(ii) Let
$\,\alpha_1<\alpha_0\,$. Then for $\,d\,$ small enough,
depending only on
$\,\alpha_1\,$ and $\,\eta\,$, the operators
$\,V(D_{\theta,\rho},u)\,$ and 
$\,b(D_{\theta,\rho},u)\,$ in $\HH$ form bounded self-adjoint analytic
families in $\theta$ provided $\,\theta\in\SS_{\alpha_1}\,$.
     \end{proposition}
{\em Proof\/}:  
 The first assertion is an immediate consequence of
Lemma~\ref{DTI properties}, Proposition~\ref{general reso p estimate}
and (a1)--(a2). For (ii)
notice that we have for a function $f$ obeying the same requirements as
$\gamma$ in (a1) and (a2)  and $\,\im\theta=\beta\,$ the bound 
\begin{equation}\label{eq:bound on DTIO}
\left\|\,f(D_{i\beta,\rho})\,\right\|\,\leq\,
\Bigl({\rm dist}
(\partial\VV,\Sigma_{|\beta|,\|\rho'\|_\infty})\Bigr)^{-1}
\int_{\partial\VV}|f(z)|\,|dz|\,,
\end{equation}
for any integration path 
$\,\partial\VV\subset\Sigma_{\alpha_0,\eta_0}\setminus\overline
\Sigma_{|\beta|,\|\rho'\|_\infty}\,$ verifying the conditions in
Definition~\ref{DTI}. 
Thus we see that the operators $\,\gamma_{\theta,\rho}\,$ 
for $\,\theta\in\SS_{\alpha_1}\,$ can be bounded by a constant depending
only on $\alpha_1$ and $\eta$. Choosing $d$ small enough, this immediately
implies that $\,\|u\gamma_{\theta,\rho}\|\,$ can be made smaller than one. 
Thus $(1+u\gamma_{\theta,\rho})^{-1}$ exists and is bounded, which is all
we need in view of the structure of $b$ and $V$ and (i).
\quad \QED
\vskip3mm
When there is no 
possibility of confusion, we will use for the operators
$\, h(D_\vartheta,u)\,$ the symbol $\, h_\vartheta\,$, too.\\
Even though formula (\ref{eq:bound on DTIO}) will be useful later on, it is not sufficient.
In particular we will need more
information about the dependence of the norm on $\,\beta\,$ which is
provided by the following proposition.
   \begin{proposition}\label{difference of the gammas}
Let $f$ obey (a1) and (a2). Then for
any compact subset $\,I\,$ of $\,(-\alpha_0,\alpha_0)\,$
there exists a constant $C$  such that for all
$\,\alpha,\beta\in I\,$,
$$
\left\|\,f(D_{i\beta})\!-\!f(D_{i\alpha})\,\right\|\,\leq\,
C \,\sin \Bigl|{\beta-\alpha\over 2}\Bigr|
\,.
$$
\end{proposition}
 {\em Proof:} We have
$$
{f_{i\beta}\!-\!f_{i\alpha} \,=\,
{i\over 2\pi}\,
\int_{\partial\VV}\,f(z)\left(r_{i\beta}(z)-r_{i\alpha}(z)\right)\,dz }\,,
$$
where $\,r_\bullet(z):= (D_\bullet-z)^{-1}$.
For $\,v,w\in L^2(\R)\,$ let $\hat v:=r _{i\beta} (z)v$
and $\hat w:= r_{-i\alpha}(\overline{z})w\,$. Then
\begin{eqnarray*}
\left((r _{i\beta} (z)-r_{i\alpha}(z))v,w\right)
&\!=\!&
(\hat v,D_{-i\alpha}\hat w)-
(D_{i\beta}\hat v,\hat w)\\ \\
&\!=\!& i\left\{\left.
\hat v\overline{\hat w}\right|^{+\omega-0}_{-\omega+0}
+ \left. e^{-i\alpha}\hat v\overline{\hat w}
\right|^{-\omega-0}_{+\omega+0}
\right\}
+\left(e^{i(\beta-{\alpha})}\!-\! 1\right)
(\chi_{\Omega_e} D_{i\beta} \hat v,\hat w)
\\ \\
&\!=\!&\left.
ie^{-i\alpha}\left(\!1\!-\!e^{-i{(\beta-\alpha)/ 2}}\right)
\hat v\overline{\hat w}\right|^{-\omega-0}_{+\omega+0}
+\left(e^{i(\beta-{\alpha})}\!-\! 1\!\right)
(\chi_{\Omega_e} D_{i\beta} \hat v,\hat w)
\,,
    \end{eqnarray*}
where we used that $\,\hat v \in \DD(D_{i\beta} )$
and that $\,\hat w \in \DD( D_{-i\alpha} )$. We can
now rewrite the boundary term
as
$\,i \hat v\overline{\hat w}|^{-\omega-0}_{+\omega+0}
=e^{i\alpha}(\chi_{\Omega_e}\hat v,D_{-i\alpha}\hat w)-
e^{i\beta}(\chi_{\Omega_e}D_{i\beta} \hat v,\hat w) \,$,
so that
$$
\bigl((r _{i\beta} (z)-r_{i\alpha}(z))v,w\bigr)\,=\,
\left(1\!-\!e^{-i{(\beta-\alpha)/ 2}}\right)
(\chi_{\Omega_e}\hat v,D_{-i\alpha}\hat w)
+
\left(e^{i{(\beta-\alpha)/ 2}}\!-\! 1\right)
(\chi_{\Omega_e}D_{i\beta} \hat v,\hat w).
$$
This can be written, dropping the argument $z$ in the resolvents, as
\begin{equation}\label{difference DTI}
{r _{i\beta}-r_{i\alpha}\over 2i} =
z\,\sin {\textstyle{\beta-\alpha\over 2}}
r_{i\alpha}\chi_{\Omega_e}r _{i\beta}+
\sin {\textstyle{\beta-\alpha\over 4}}
\left(e^{-i{(\beta-\alpha)/ 4}}
\chi_{\Omega_e}r _{i\beta}+
e^{i{(\beta-\alpha)/ 4}} r_{i\alpha}\chi_{\Omega_e}
\right)
\end{equation}
implying by proposition~\ref{general reso p estimate}
\begin{eqnarray*}
\left\|r _{i\beta} (z)-r_{i\alpha}(z)\right\|
&\!\leq \!& 2\left|\sin  {\textstyle{\beta-\alpha\over 2}}\right|
\Bigl(\|r _{i\beta} (z)\|+\left\|r_{i\alpha}(z)\right\|+
|z|\,\|r _{i\beta} (z)\|\,\left\|r_{i\alpha}(z)\right\|
\Bigr)
\\ \\&\!\leq \!& 2\left|\sin {\textstyle {\beta-\alpha\over 2}}\right|
{\mbox{dist}(z,\Sigma_{|\beta|,0})+|z|+\mbox{dist}(z,\Sigma_{|\alpha|,0})\over
\mbox{dist}(z,\Sigma_{|\beta|,0})\,\mbox{dist}(z,\Sigma_{|\alpha|,0})}
\,.
\end{eqnarray*}
We also used that since $|\beta-\alpha|<\pi$, we can bound $\vert
\sin{\beta-\alpha\over4}\vert$ by $\vert
\sin{\beta-\alpha\over2}\vert$. Furthermore we can suppose without
restriction of generality that
$|\beta|\geq|\alpha|$. Then
we can choose $\,\partial\VV\,$ in
$\,\Sigma_{\alpha_0,\eta_0/2}\setminus\Sigma_{|\beta|,0}\,$
such that the last factor on the left side is
bounded by  a constant depending only on $\,I\,$;
due to (a2) the statement follows.
\quad \QED

\subsection{Some estimates on $\,p_{i\beta}\,$ and $\,W_{i\beta,\rho}\,$}
The leading longitudinal part of the dilated Hamiltonian is the
operator of multiplication by $\,p^2_{i\beta}\,$. Here we collect
some simple bounds which we shall need in the following.

      \begin{proposition} \label{p_theta estimates}
(i) For all $p\in\R$, any positive integer $\,n\,$, any
$\;|\beta|\leq {\pi\over 2}\,$ and $\,\omega\!\geq\!0\,$
we have
$\,
p^{2n}\geq |p^n_{i\beta}|^2\ge  p^{2n}\cos^n\beta\,.
$
(ii) The function $\,p\mapsto p_{i\beta}(p)\,$, satisfies on
$\,\Omega_e\,$ for any $\,\omega\!\geq\!0\,$ the bounds:
  $\;\re p_{i\beta}^2\ge p^2\cos 2\beta+
2\omega^2\sin^2\beta\;$ if $\;|\beta|\,\leq\, {2\pi\over 3}\,$.
    \end{proposition}

\noindent
{\em Proof:} (i) It is sufficient to consider $\,n=1\,$.
For every real $p$ we have
$$
|p_{i\beta}|^2\,=\, p^2\chi_{\Omega_i}\,+\,
\left(\omega^2+(|p|\!-\!\omega)^2+ 2\omega(|p|\!-\!\omega) \cos\beta
\right)\chi_{\Omega_e}\,.
$$
The part on the \rhs restricted to $\,\Omega_e\,$ satisfies
$$
p^2+2(\omega^2-\!\omega |p|)(1\!-\!\cos\beta) \,\ge\,
p^2+(\omega^2-\!p^2)(1\!-\!\cos\beta) \,=\, p^2\cos\beta+\omega^2
(1\!-\!\cos\beta)\;;
$$
since the very last term is nonnegative,
we obtain the second inequality.
Furthermore, $\,\omega^2-\!\omega |p|<0\,$ on
$\,\Omega_e\,$, so the same expression may be estimated from above
by $\,p^2\,$ and thus the first inequality follows.

The identity
$$
\re p_{i\beta}^2\,=\, p^2\cos 2\beta+ \omega^2 (1\!-\!\cos 2\beta)+
2\omega (|p|\!-\!\omega)(\cos\beta\!-\! \cos 2\beta)
$$
yields (ii) on $\,\Omega_e\,$, because $\,\cos\beta\!-\! \cos 2\beta
\ge 0\,$ for $\;|\beta|\le\, {2\pi\over 3}\,$.  \quad \QED
\vspace{3mm}

\noindent
Let us fix an $\alpha_1$ in $\,(0,\alpha_0)\,$ and define the weight
\begin{equation}\label{good weight}
\,\langle p \rangle \,:=\,
\left(p^2 +\tau \right)^{1/2}\,,\quad\tau\,:=\,
\sup\{\|V^0_{i\beta,\rho}\|\,:\,|\beta|\leq\alpha_1\,,\
\rho\mbox{\ verifying (\ref{condition on rho})}\}\,,
\end{equation}
The supremum exists by (\ref{eq:bound on DTIO})
and is strictly positive by (a0).

The motivation for the choice of this weight is that it will simplify
the statements and permit us to obtain
particularly simple
constants in the subsequent resolvent estimates. Notice that it
depends only on the fixed parameters $\,\alpha_1\,$ and $\,\eta\,$, but
not on $\,d\,$.

Furthermore we fix $d_0$ such that $b_{i\beta,\rho}$ and 
$V_{i\beta,\rho}$ exist and are bounded for all $|\beta|\leq\alpha_1$ and 
all $\rho$ verifying (\ref{condition on rho}), if $d\leq d_0$.
     \begin{proposition}\label{W_thetarho estimate}
Let $\,|\beta|<\alpha_1\,$, and  $\,d \leq d_0\,$.
Then there exists a  constant $\,c_{i\beta}^{W}>0\,$
for all $\,\rho\,$ satisfying (\ref{condition on rho})
such that $\,
\left\|\langle p \rangle^{-1} W_{i\beta,\rho}\langle p \rangle^{-1}\right\|
\le c_{i\beta}^{W}d $.
   \end{proposition}
{\em Proof:} Since
$\,
W _{i\beta,\rho} \,= \,p_{i\beta}
(b\!-\!1)_{i\beta,\rho}p_{i\beta}
+ (V\!-\!V^0)_{i\beta,\rho}
\,$
we get 
$$
c_{i\beta}^{W} =\max_{0\leq d\leq d_0}\left\{ {1\over d}\left(
\|(b\!-\!1)_{i\beta,\rho}\| +
 {1\over\tau}\|(V\!-\!V^0)_{i\beta,\rho}\|\right)\right\}
$$
which does exist because
$\,|(b\!-\!1)(\cdot\,,u)|d^{-1}\,$ and
$\,|(V\!-\!V^0)(\cdot\,,u)|d^{-1}\,$
obey (a2) for all $u\in[0,d],\;d\in(0,d_0]$.
Notice that we also used $\,|p_{i\beta}/p|\leq 1\,$ as proven
in Proposition~\ref{p_theta estimates}.
 \quad \QED
\subsection{The operators $\,H^0_{\theta}\,$ and $\,H_{\theta}\,$}
Let us finally collect some basic properties of the operators
 $\,H^0_{\theta}\,$ and $\,H_{\theta}\,$, images under
 the exterior scaling of the ``free" and the full Hamiltonian,
 $\,H^0\,$ and $\,H\,$, respectively.
      \begin{theorem} \label{scaled Hamiltonians}
(i) The operators $\,H^0_{\theta}\,$ for
$\,|\im\theta|<\alpha_0\,$ form a self--adjoint analytic family
of type A with the common domain
$\,\DD(H^0)=\DD(p^2)\otimes (\HH_0^{1}\cap \HH^{2})((0,d))$.
Moreover (cf. Figure~3),

$$
\sigma(H^0_{\theta})\,=\, \left\lbrace\, \lambda+E_j\,:\;
\lambda\in \left(\{\lambda_n\}_{n=1}^N \cup\nu\cup \sigma(p^2_{\theta})
\right)\,, \; j=1,2,\mbox{\dots}\; \right\rbrace\,,
$$
where $\,\nu\,$ denotes the set of resonances of the operators
$\,A_\theta =p^2_\theta + V^0_\theta\,$ (which may be empty).

\noindent
(ii) For all sufficiently small
$\,d\,$,
the operators $\,H_{\theta}\,$ with $\,|\im\theta|<\alpha_1\,$
 form a self-adjoint analytic family of type A
with the common domain $\,\DD(H^0)\,$.
      \end{theorem}
      \begin{figure}

\epsfig{file=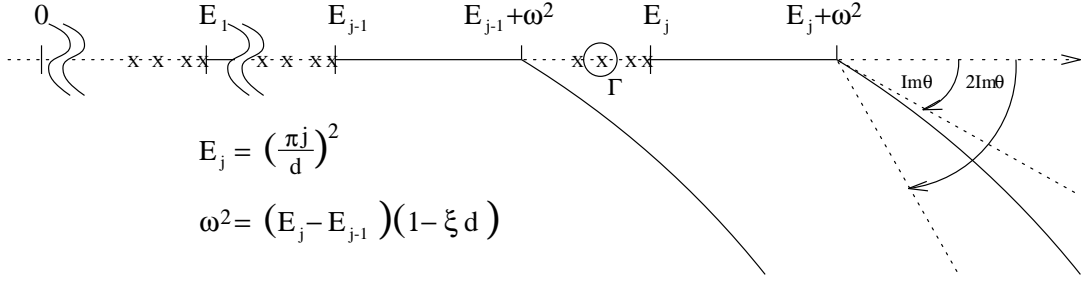}
   \caption{The spectrum of $\,H^0_{\theta}\,$ ---
Stability of the eigenvalues}
      \end{figure}
\noindent
{\em Proof:}
(i) Analyticity of $\,H^0_{\theta}\,$ follows from the boundedness
 and the analyticity of $\,\gamma_{\theta}\,$, see
Proposition~\ref{DTI gamma}. The form of the  spectrum
 is  due to the analyticity, the structure of the operator,
$\,H^0_{\theta} \, =\, A_\theta\otimes I + I\otimes (-\partial_u^2)\,$,
and the $\,p^2$-compactness of $\,\gamma_{\theta}^2$
in view of (a2).\\
(ii) Due to Proposition~\ref{W_thetarho estimate}, $\,W_{i\beta}$
is $\,H^0_{i\beta}$--bounded with a bound smaller than one if
$\,d\,$ is small enough and $\,|\beta|<\alpha_1\,$. This bound extends
by the group property of 
the exterior dilation to all $\,\theta\in\SS_{\alpha_1}\,$. Thus the
analyticity of  
$\,H_{\theta}\,$ follows.
 \quad \QED
\begin{remark}{\rm 
(i) The set $\,\nu\,$ of resonances of $\,A_\theta\,$
cannot contain embedded
eigenvalues due to the decay assumption (a2)
(\cf {\rm \cite[Sect. XIII.13]{RS}}).\\
(ii) If $\,\omega\,$ is chosen large enough no resonances
of $\,A_\theta\,$ will be disclosed at all;
the continuous spectrum is deformed only in a sector with vertex
$\,\omega^2\,$, whereas the resonances lie inside a disc around
the origin with a fixed radius of order $\,\|V^0\|\,$.}
\end{remark}
\setcounter{equation}{0}
\section{Stability of the discrete spectrum}           
Our next goal is to estimate the effect the perturbation
$\,W_{i\beta}\,$ has on the spectrum of the operator
$\,H^0_{i\beta}\,$.
\subsection{Estimates on $\,R^0_{i\beta}\,$ and $\,R^0_{i\beta,\rho}\,$}

Let $\,E^0_{j,n}=\lambda_n+E_j\,,j\geq 2\,$
 be a fixed eigenvalue of $\,H^0_{\theta}\,$.
We choose:
\begin{equation}\label{fixed Gamma}
\,\Gamma:=\{z\in\C:\: |z\!-\!E^0_{j,n}|=r\,\},\;
\mbox{ with}\;
r:={1\over 2}{\rm dist}
(\lambda_n,\sigma(A)\setminus\{\lambda_n\})
\end{equation}
to be
a circular contour
around $\,E^0_{j,n}\,$
such that no other eigenvalue of $\,H^0_{i\beta}\,$
is within $\,\Gamma\,$, and denote $\,\DD_{\Gamma}:=\{z\in\C:\:
|z\!-\!E^0_{j,n}|\le r\,\}\,$.

Having the intention to prove stability by perturbation we have to
control $\,R^0_{i\beta}(z)\,$ on $\,\Gamma\,$.
For the estimate it is advantageous to pass to the
transverse mode decomposition, $\,H^0_{i\beta}\,$ being
diagonal in this decomposition:
$$
H^{0}_{i\beta}=\sum_{k\geq 1}\JJ_kH^{0,k}_{i\beta}\JJ_k^*,\qquad
 H^{0,k}_{i\beta}:= \JJ_k^* H^0_{i\beta} \JJ_k
\;\;{\rm on}\;\; L^2(\R,dp)\,.
$$
Since $\, H^{0,k}_{i\beta}\,$ is not self-adjoint (for $\beta\neq 0$)
and a part of $\Gamma$ lies in the numerical range of
$\, H^{0,k}_{i\beta}\,$ for every $k\leq j\:$ (\cf Figure 3), 
simple estimates in terms of the distance to the numerical range
do not work here.

The difficulties for the  $j$-th mode result from our desire to chose
the $\,\omega\,$ of the exterior dilation (see eq.(\ref{dilation}))
as $\,\omega=\OO(\sqrt{E_j-E_{j-1}})\,$.
This $d$ dependence of $\,\omega\,$ implies a
$d$ dependence of the deformed operators $\,H^{0,k}_{i\beta}\,$ so that
the usual argument using the compactness of $\Gamma$ to assure uniform
boundedness of $\,R_{i\beta}^{0,j}(z)\,$ are not applicable here.
Instead we chose to perturbate around $\,\beta=0\,$; it turns out
that then the estimate is independent of $\,\omega\,$, and thus of
$\,d\,$: 
\begin{lemma}\label{lem reso Rj estimate}
Let 
$\,c^{(2)}:={2\over r} \max\{1,3\tau\}\,$.
There exists
$\,0<\beta_1\leq\min\{{\pi\over 4},\alpha_1\}\,$
so that for $\,|\beta|\leq \beta_1\,$
one has
$$\;\sup_{z\in\Gamma} \|\langle p \rangle^\ell R_{i\beta}^{0,j}(z)
\langle p \rangle\|
\le c^{(2)},\,l=0,1\quad {\rm and}\quad
\sup_{z\in\Gamma} \|\langle p \rangle^\ell R_{i\beta}^{0,j}(z)\|
\le c^{(2)},\,l=0,1,2.
$$
\end{lemma}
{\em Proof:\/}
We prove the claim using that $\,R_{i\beta}^{0,j}(z)\to
R^{0,j}(z)\,$ in the operator norm as $\,\beta\to 0\,$ uniformly for
$\,z\in\Gamma\,$. By the resolvent identity,
$$
R_{i\beta}^{0,j}(z)-R^{0,j}(z)\,=\, R_{i\beta}^{0,j}(z)
\left(p^2\!-p^2_{i\beta}+V^0\!-V^0_{i\beta}\right)R^{0,j}(z)\,.
$$
One has for every $\,\beta\,$ the inequality
$|p^2\!-p^2_{i\beta}|\,\leq 6\, |\sin{\beta\over 2}|p^2\,$.
Let $\,0<\beta_1<\alpha_0\,$.
Using Proposition~\ref{difference of the gammas} implies that
there is a constant $\,C\,$ for all
$\,|\beta|\leq\beta_1\,$  such that
$\,\|V^0\!-\!V^0_{i\beta}\| \,\le\, C\, |\sin{\beta\over 2}|\,$. Thus
we get the
estimate for every $\,z\in\Gamma\,$:
      \begin{eqnarray*}
\lefteqn{\|\left(p^2\!-p^2_{i\beta}+V^0\!-V^0_{i\beta}\right)R^{0,j}(z)\|}
\\ \\ &&\,=\,\left\|\,\left\lbrace\, \left(p^2\!-p^2_{i\beta}\right)
(p^2\!+\!1)^{-1} +(V^0\!-V^0_{i\beta}) (p^2\!+\!1)^{-1} \right\rbrace\,
(p^2\!+\!1)\,R^{0,j}(z)\, \right\| \\ \\ &&
\,\le\, \max\{6,C\}
|\sin{\beta\over 2}|
\, \left(\,
1+\left( 1\!+\!r\!+\!\|V^0\!-\!\lambda_n\|\right)
{1\over r}
 \right)\,.
      \end{eqnarray*}
Taking $\,\beta_1\,$ small enough one can certainly ensure that
for all $\,\beta\,$ with $\,|\beta |\leq\beta_1\,$ and all
$\, z\in \Gamma\!$ one has
$\|\left(p^2\!-p^2_{i\beta}+V^0\!-V^0_{i\beta}\right)R^{0,j}(z)\|
\leq{1\over 2}$. Solving the resolvent identity gives
$\,\|R^{0,j} _{i\beta}(z)\|\leq {2\over r}\,$.
Thus the use of
$\,\langle p \rangle^2R^{0,j}(z)=1+\left(\tau\!+\!z\!-\!E_j\!-\!V^0 \right)
R^{0,j}(z)\,$, and of the facts that in our situation
$\|V^0-\lambda_n\|\leq \|V^0\|$ and $r\leq \|V^0\|/2$,
yields the claim for $\ell=2$ in the second formula.
The statement in the case of only one weight present is
then obtained by an obvious quadratic estimate.
The symmetric case with one weight on each side of the resolvent
is handled by following estimate:
$$
\|\langle p \rangle R^{0,j}(z)\langle p \rangle\|^2\! \leq\!
\|\langle p \rangle R^{0,j}(z)\langle p \rangle\| +
(r+2\tau)\|\langle p \rangle R^{0,j}(z)\|^2.
\quad\QED
$$
The restriction $\,|\beta|<{\pi\over 4}\,$ is not necessary in this
lemma, but for later convenience we prefer having it stated.
Indeed to simplify the statements we work from
now on with a {\em fixed\/} $\beta$:
\begin{equation}\label{fixed beta}
\beta\in[-\beta_1,0),\;{\rm with}\;\beta_1\; \mbox{given by Lemma~4.1}.
\end{equation}
Consequently, the dependence of the constants on $\beta$ will be
no longer specified. 

For $\,k\neq j\,$ the resolvents $\,R^{0,k}_{i\beta}(z)\,$ are
estimated considering $\,V^0_{i\beta}\,$ as a perturbation.
We use the fact that
the distance between $\,\DD_\Gamma\,$
and the spectrum of $\,H^{00,k}_{i\beta}:=E_k+p^2_{i\beta}\,$ tends
to infinity as $d$ tends to zero whereas
$\,V^0_{i\beta}\,$ is bounded, independently of $d$.
We choose $\,\omega\,$
to be
\begin{equation} \label{dilation threshold}
\omega\,:=\, {\pi\over d}\, \sqrt{(2j\!-\!1)(1\!-\!\xi d)}\,,
      \end{equation}
where
$\,\xi\,$ is a supplementary positive parameter specified below to
govern the distance of the spectrum of
$\,H^{00,j-1}_{i\beta}\,$ to the contour $\,\Gamma\,$ --- \cf Figure 3.
 \begin{lemma}\label{resolvent lemma}
Let $\rho$ verify condition (\ref{condition on rho})
 and
$\,c^{(1)}:=8\sqrt{3}\,$.
Then
for all $\,(d,\xi)\,$ verifying
\begin{equation}\label{second domain xi d}
\,1\,\geq\,\xi d \,
\geq
\,{c^{(1)}\,\tau\over |\sin\beta|}\,d^2\,
\end{equation}
one has for $\,\ell =0,1\,$
                \begin{description}
   \item{(a)} $\;\sup_{z\in \DD_{\Gamma}} \|\langle p \rangle^\ell
R_{i\beta,\rho}^{0,j-1}(z)\langle p \rangle\|
\le {c^{(1)}\over |\sin\beta|\xi}d^{-\ell}$,
   \vspace{-.2em}
   \item{(b)} $\;\sup_{z\in \DD_{\Gamma}} \|\langle p \rangle^\ell
R_{i\beta,\rho}^{0,k}(z)\langle p \rangle\|
\le {c^{(1)}\over |\sin\beta|}d^{1-\ell}\,$
for all $\,k\ne j,\, j\!-\!1\,$.
                \end{description}
   \vspace{-.2em}
\end{lemma}
{\em Proof of the lemma} is given in appendix B.
\subsection{Stability of the resolvent set of $\,H^0_{i\beta}\,$}
   \begin{lemma} \label{resolvent set stability}
Let $\,\Gamma\,$ and
$\,\beta\,$ be fixed by (\ref{fixed Gamma}) and (\ref{fixed beta}).
Then for all sufficiently small $\,d\leq d_0\,$
such that the condition (\ref{second domain xi d}) is verified for
$\,\xi\geq 2c^W
\max\{c^{(2)},{c^{(1)}\over|\sin{\beta}|}\}
\,$ ($c^{W}\equiv c^{W}_{i\beta}$), the contour
$\,\Gamma\,$ belongs to the resolvent set $\,\rho(H_{i\beta})\,$.
   \end{lemma}

\noindent
{\em Proof:} If we can show that
$$
R_{i\beta}(z)\,=\, R^0_{i\beta}(z)\langle p \rangle\left
(1\!+\!\langle p \rangle^{-1} W_{i\beta}R^0_{i\beta}(z)\langle p \rangle
\right)^{-1}\langle p \rangle^{-1}
$$
makes sense for all $\,z\in\Gamma\,$, we are done. The boundedness of
$\,R^0_{i\beta}(z)\langle p \rangle\,$ has  already been proven
in the two preceding lemmas. Thus showing that
$\,\|\langle p \rangle^{-1}W_{i\beta}
R^0_{i\beta}(z)\langle p \rangle \| <1\,$ will conclude the proof.
First of all one has $\,R^0_{i\beta}(z)= \sum_{k\geq 1} \JJ_k
R^{0,k}_{i\beta,\rho}(E) \JJ^*_k$.
Secondly, the operators $\,\JJ_k,\, \JJ_k^*\,$
commute with $\,\langle p \rangle\,$ and the map
$\,\JJ:=\sum_{k\geq 1} \JJ_k:\: \bigoplus_{k\geq 1} L^2(\R,dp) \to\HH\,$
is an isometry. So we can employ Lemmas~\ref{lem reso Rj estimate},
\ref{resolvent lemma}
and Proposition~\ref{W_thetarho estimate} to get
      \begin{eqnarray*}
\|\langle p \rangle^{-1}W_{i\beta}R^0_{i\beta}(z)\langle p \rangle\|
\!&\le &\!
\|\langle p \rangle^{-1}W_{i\beta}\langle p \rangle^{-1}\|\,
\left\|\, \bigoplus_{k=1}^{\infty}\,
\langle p \rangle R^{0,k}_{i\beta}(z)\langle p \rangle \right\| \\ \\
\!&\le &\! c^W d\: \max_k \left\|\,
\langle p \rangle R^{0,k}_{i\beta}(z) \langle p \rangle\right\| \\ \\
\!&= &\! c^W d\:
\max \left\lbrace c^{(2)},
\,{c^{(1)}\over |\sin\beta| \xi d}\, \right\rbrace 
\,\leq\,{1\over 2}\;;
      \end{eqnarray*}
recall that $\,0<\xi d\leq 1\,$. \quad \QED

   \begin{corollary} \label{stability corollary}
Under the same assumptions, the eigenvalue $\,E^0_{j,n}\,$
of $\,H_{i\beta}^0\,$
gives rise to a single perturbed eigenvalue of $\,H_{i\beta}\,$ of
the same multiplicity.
   \end{corollary}

\noindent
{\em Proof:} By standard interpolation between the respective projections,
$$
P^0_{i\beta}\,:=\,{i\over 2\pi}\, \int_{\Gamma}\, R^0_{i\beta}(z)\,dz
\quad {\rm and} \quad P_{i\beta}\,:=\,{i\over 2\pi}\, \int_{\Gamma}\,
R^0_{i\beta}(z) (1\!+\!W_{i\beta}R^0_{i\beta}(z))\,dz\,. \qquad \QED
$$
\setcounter{equation}{0}
\section{Exponential decay estimates}
Please keep in mind that
$\,\beta\,$ is now considered to be a fixed parameter, \cf
(\ref{fixed beta}) and that, up to the end of the proof of Theorem~2.2, 
$\rho$  obeys the condition (\ref{condition on
rho}); these facts might not always be stated explicitely.
Let $\,E\equiv E_{j,n}\,$
be the resonance associated with $\,E_{j,n}^0\,$.
Under the conditions of the
last section
 the corresponding complex eigenvalue equation
\begin{equation}\label{eq: complex ev}
H_{i\beta}\phi_{i\beta}\,=\, E\phi_{i\beta}
\end{equation}
can be easily demonstrated to be equivalent to the system
                        \begin{eqnarray}\label{Feshbach 1}
\left( P_jH_{i\beta}P_j- P_jW_{i\beta}\widehat{R}^j_{i\beta}(E)
W_{i\beta}P_j \right)\phi_{i\beta} &\!=\!& EP_j\phi_{i\beta}\,, \\ \nonumber \\
\label{Feshbach 2}
Q_j\phi_{i\beta} &\!=\!& -\widehat{R}^j_{i\beta}(E) W_{i\beta}P_j\phi_{i\beta}
                        \end{eqnarray}
for a given $\,j=2,3\dots\,$, as pointed out in Section~2.3. Recall
that there we defined $\,\widehat{R}^j_{i\beta}(E)=
Q_j(Q_j(H_{i\beta}\!-\!E)Q_j)^{-1}Q_j\,$.
We shall introduce the analogous notation, $\,\widehat{A}^j\,$, also 
for a general closed operator $\,A\,$: we define $\,\widehat{A}^j:=
Q_jAQ_j\,$ meaning that the operator is restricted to the orthogonal
complement of the subspace associated to the mode $\,H^{0,j}\,$. In
the case of resolvents the hat designates the resolvent on $\,Q_j\HH\,$,
that is $\,\widehat{(A\!-\!z)^{-1}}:=Q_j(Q_j(A-z)Q_j)^{-1}Q_j\,$.

        Using the embedding operators (see Section~2.1),
we find that (\ref{Feshbach 1}) is further equivalent to
   \begin{equation} \label{eigenvalue equation}
\left( H^j_{i\beta}\!-\!B^j_{i\beta}(E) \right)\phi^j_{i\beta}
\,=\, E\phi^j_{i\beta}\,,\quad
B^j_{i\beta}(E)\,:=\,  
\JJ_j^*W_{i\beta} \widehat{R}^j_{i\beta}(E)W_{i\beta}\JJ_j 
   \end{equation}
on $\,L^2(\R)\,$. First we have to establish that these equations
make indeed sense.
      \begin{proposition} \label{function B}
Under the conditions of Lemma~\ref{resolvent set stability}
on $\,d\,$ and $\,\xi\,$ and the condition
(\ref{condition on rho}) on $\,\rho\,$ the following bounds are valid\\
(i)
$\,\|\langle p \rangle \widehat{R}^j_{i\beta,\rho}(E)\langle p \rangle \|
 \le {2c^{(1)}\over |\sin\beta|\xi d}\,$,\\
(ii)
$\,\|\langle p \rangle^{-1}W_{i\beta,\rho}\langle p \rangle^{-1}\|\,
\|\langle p \rangle\widehat{R}^j_{i\beta,\rho}(E)\langle p \rangle\|
\leq 1\,
$
and,\\
 (iii) $\,\|\langle p \rangle^{-1}B^j_{i\beta,\rho}(E)\langle
p \rangle^{-1}\|
 \le c^W d\,$.
      \end{proposition}
\noindent
{\em Proof:} We can write
\begin{equation}\label{reduced reso eq}
\widehat{R}^{j}_{i\beta,\rho}(E)
\,=\, \left(Q_j H^0_{i\beta,\rho}Q_j\!+\!Q_j
W_{i\beta,\rho}Q_j\!-\!E\right)^{-1}
\,=\, \widehat{R}^{0,j}_{i\beta,\rho}(E)
\left(I\!+\!\widehat{W}_{i\beta,\rho}
\widehat{R}^{0,j}_{i\beta,\rho}(E) \right)^{-1}Q_j\,.
\end{equation}
Now since
$\,\widehat{R}^{0,j}_{i\beta,\rho}(E)= \sum_{k\ne j} \JJ_k
R^{0,k}_{i\beta,\rho}(E) \JJ^*_k\,$
one has by the argument in the proof of
Lem\-ma~\ref{resolvent set stability} and by Lemma~\ref{resolvent lemma}
  \begin{eqnarray*}
\left\|\langle p \rangle^{-1}\widehat{W}_{i\beta,\rho}
\widehat{R}^{0,j}_{i\beta,\rho}(E)\langle p \rangle\right\| \!&\le &\!
c^W d\: \max_{k \ne j}
\left\|\langle p \rangle R^{0,k}_{i\beta,\rho}(E)\langle p \rangle\right\|
\,\leq\, \frac{1}{2}\,.
   \end{eqnarray*}
Hence
$\,\|\langle p \rangle\widehat{R}^j_{i\beta,\rho}(E)\langle p \rangle\|
\leq {2c^{(1)}\over |\sin\beta|\xi d}\,$.
The condition on $\,\xi\,$ ensures then that
$\,{2c^{(1)}c^W\over |\sin\beta|\xi}\leq 1$
and thus (ii) by
Proposition~\ref{W_thetarho estimate}.
The last assertion is due to the estimate
$$
\left\|\langle p \rangle^{-1}B^j_{i\beta,\rho}(E)\langle p \rangle^{-1}
\right\|\leq
\left\|\langle p \rangle^{-1}W_{i\beta,\rho}\langle p \rangle^{-1}\right\|^2
\left\|\langle p \rangle\widehat{R}^j_{i\beta,\rho}(E)\langle p \rangle
\right\|
\leq\,c^W d.\quad\QED
$$
\noindent
Let $\,\phi_{i\beta}\,$ be a normalized solution of the above complex
eigenvalue equation (\ref{eq: complex ev}).
Denote the boosted eigenfunction $\,e^\rho\phi_{i\beta}\,$
by  $\,\phi_{i\beta,\rho}\,$, where $\,\rho$ 
obeys (\ref{condition on rho}).
Then equation (\ref{eigenvalue equation}) implies
$$
e^\rho (H^j_{i\beta} -B^j_{i\beta} (E)-E)e^{-\rho}\phi^j_{i\beta,\rho}\,=\,0\,,
$$
which in turn gives the relation
\begin{equation}\label{fundamental eq for decay estimate}
\re\left(\left(H^j_{i\beta,\rho} -B^j_{i\beta,\rho} (E)-E\right)
\phi^j_{i\beta,\rho},\phi^j_{i\beta,\rho} \right)\,=\,0.
\end{equation}
To be able to apply the usual Agmon technique \cite{Ag},
we need the following
\begin{proposition}\label{form sense estimate for Agmon}
Let $\,d\,$ be small enough so that $\,\cos 2\beta -2c^Wd>0\,$.
Then $\,p_\star\,$ defined by
\begin{equation}\label{sStar}
{1\over 2}p^2_\star \,:={\|V^0_{i\beta,\rho}\|\!
 +\! |\re E\!-\!E_{j}|\!+\!2c^W \,\tau\,d
\over
(\cos 2\beta -2c^Wd)},
\end{equation}
 is uniformly bounded in $\, d\,$.
Under the conditions of Proposition~\ref{function B}
the following inequality holds in the form sense on $\DD(p^2\!\otimes\! I)$:
$$
\re\left(\,H^j_{i\beta,\rho} -B^j_{i\beta,\rho} (E)-E\right)
\,\geq\,
(\cos 2\beta -2c^Wd)\left(
p^2-{1\over 2}p^2_\star\right)\,.
$$
\end{proposition}
{\em Proof:\/}
The statement on $\,p_\star\,$ is trivial. Then
using the estimates on $\,p_{i\beta} ,W_{i\beta,\rho}\,$ and
$\,B^j_{i\beta,\rho}$ of Propositions \ref{p_theta estimates},
\ref{W_thetarho estimate}, and \ref{function B}
we get
   \begin{eqnarray*}
\re\Bigl(H^j_{i\beta,\rho}\!-\!B^j_{i\beta,\rho} (E)\!-\!E\Bigr)
&\!\geq\!& \re p^2_{i\beta} +
\re\!\left(W^j_{i\beta,\rho}\!-\!B^j_{i\beta,\rho}(E)\right)
-\|V^0_{i\beta,\rho}\|\!
-\re (E \!-\!E_j)\,
 \\ \\
&\!\geq\!& \cos(2\beta)p^2 \! -2c^Wd\langle p \rangle^2
- \|V^0_{i\beta,\rho}\|\! -|\re E\! -\!E_{j}|\,.
\quad\QED
   \end{eqnarray*}
We conclude this section with the main
ingredient for the proof of the estimate on the resonance width:
the exponential decay as $\, d\,$ tends to zero of the resonance wave
function $\, \phi^j_{i\beta}\,$ in the $\HH^1$ sense.
  \begin{theorem} \label{decay estimates thm}
 Denote
$\,\Omega_\star:= (-p_\star,p_\star)\,$, where $\,p_\star\,$ is
defined by (\ref{sStar}).
Then  for any $\,\beta\in[-\beta_1,0)\,$
and any $\,\eta\in (0,\eta_0)\,$
there is a $\,d_\eta\leq d_0\,$,
 such that for $\, d\in(0,d_\eta)\,$ one has
$\,\cos 2\beta -2c^Wd>0\,$ and $\,\omega>p_\star\,$,
with $\,\xi\,$ as in Lemma~\ref{resolvent set stability},
and for
\begin{equation}\label{eq:def of rho}
\rho(p)\,:=\, \eta\, \int_{\min\{0,p\}}^{\max\{0,p\}}
\chi_{\Omega_i\setminus \Omega_\star}(t)\, dt.
\end{equation}
we have
\begin{equation}\label{eq:reso fct estimate}
\left\|\phi^j_{i\beta,\rho}\right\|^2\,\leq\, 2 \qquad {\rm and} \qquad
\left\|p\phi^j_{i\beta,\rho}\right\|^2\,\leq\, 2p_\star^2 \,.
\end{equation}
   \end{theorem}
{\em Proof:\/}
The first statement is evident, since $\beta$ and $\eta$ are fixed
parameters and $p_\star$ remains bounded as $d$ tends to zero.
For the proof of the second statement note
that $\,\rho\,$ satisfies (\ref{condition on rho}). At the same time,
$\,\rho'\,$ is by 
definition zero on $\,\Omega_\star\subset\Omega_i\,$.
So
 the use of the preceding proposition  with
$\,\cos 2\beta -2c^Wd>0\,$ and
the relation (\ref{fundamental eq for decay estimate}) yields
\begin{eqnarray*}
\left(\left(p^2 \!-\!{1\over 2}p^2_\star\right)\chi_{\Omega_\star^c}
\phi^j_{i\beta,\rho},\phi^j_{i\beta,\rho}\right)
&\,\leq \,&
\left(\left({1\over 2}p^2_\star\!-\!p^2\right)\, \chi_{\Omega_\star}
\phi^j_{i\beta,\rho},\phi^j_{i\beta,\rho}\right)\\ \\
&\,\leq \,& {1\over 2}p^2_\star \left\| \chi_{\Omega_\star}
\phi^j_{i\beta}\right\|^2\leq{1\over 2}p^2_\star
\end{eqnarray*}
Evidently we have
$$
(p^2 \!-\!{1\over 2}p^2_\star)\,\chi_{\Omega_\star^c}
\,\ge\,{1\over 2}p^2_\star\,\chi_{\Omega^c_\star}.
$$
Inserting this into the above inequality, we first find
$\,\|\chi_{\Omega_\star^c}\phi^j_{i\beta,\rho}\|^2\le 1\,$, and
using the same inequality for the second time, we arrive at the
estimate
$$
\|p\chi_{\Omega_\star^c}\phi^j_{i\beta,\rho}\|^2\,\le\,
p^2_\star\,.
$$
The observation that
$\,\left\|\phi^j_{i\beta,\rho}\chi_{\Omega_\star}\right\|=
\left\|\phi^j_{i\beta}\chi_{\Omega_\star}\right\|\leq 1\,$ and that
$\,\left\|p\phi^j_{i\beta,\rho}\chi_{\Omega_\star}\right\|\leq p_\star\,$
finishes the proof.\quad\QED

\setcounter{equation}{0}
\section{Concluding the proof of Theorem~2.2}

Up to now we have employed the real part of equation (\ref{eigenvalue
equation}). The imaginary part yields
\begin{equation}\label{eq 6.1}
\im E\|\phi^{j}_{i\beta}\|^2 =
(\im (H^j_{i\beta} -B^j_{i\beta}(E))\phi^{j}_{i\beta},\phi^{j}_{i\beta})\;;
\end{equation}
for the moment we do not need the complex boosts. Using the following
simple identity,
   \begin{equation}\label{im ABA}
        \im (ABA) = 2\re\left[\,\im (A) BA\,\right]+ A^*\,\im (B) A\,,
   \end{equation}
together with resolvent equation, we can express $\,\im B^j_{i\beta}\,$,
 as
                \begin{eqnarray*}
\im B^j_{i\beta} &\!=\!& Z_{i\beta} +
\im E\:|\widehat{R}^j_{i\beta}{W}_{i\beta}\JJ_j|^2, \\ \\
Z_{i\beta} &\!:=\!& \JJ^*_j\left\{
2\re\left[\im (W_{i\beta})\widehat{R}^j_{i\beta}{W}_{i\beta}\right]
-{W}^*_{i\beta}\widehat{R}^{j\;*}_{i\beta}
\im (\widehat{H}^j_{i\beta})\widehat{R}^j_{i\beta}{W}_{i\beta}
\right\}\JJ_j,
        \end{eqnarray*}
where we have already ceased denoting the explicit dependence of
the resolvents on $\,E\,$.
Inserting this into (\ref{eq 6.1}) we get
$$
\im E\left(\|\phi^{j}_{i\beta}\|^2 +
\|\widehat{R}^j_{i\beta}{W}_{i\beta}\JJ_j\phi^{j}_{i\beta}\|^2 \right)\,=\,
((\im H^j_{i\beta} \!-\!Z_{i\beta})\phi^{j}_{i\beta},\phi^{j}_{i\beta})\,.
$$
Now the equation (\ref{Feshbach 2}) together with $\,\JJ_j \JJ^*_j
= I_{L^2(\R,dp)}$ yields
$$
\|\phi^{j}_{i\beta}\|^2 +
\|\widehat{R}^j_{i\beta}{W}_{i\beta}\JJ_j\phi^{j}_{i\beta}\|^2 \,=\,
\|P_j\phi_{i\beta}\|_\HH^2+
\|Q_j\phi_{i\beta}\|^2_\HH \,=\, \|\phi_{i\beta}\|^2_\HH\;;
$$
hence if the complex-scaled eigenvector $\,\phi_{i\beta}\,$ is normalized,
equation (\ref{eq 6.1}) is equivalent to
\begin{equation}\label{fundamental eq for ImE}
\im E\,=\, ((\im H^j_{i\beta} \!-\!Z_{i\beta})\phi^{j}_{i\beta},
\phi^{j}_{i\beta})\,.
\end{equation}
The following proposition shows that considering the imaginary part
means in a sense a localization of the Dunford-Taylor operators on
$\,\Omega_e\,$, which is naturally the case for local operators, \ie
\hskip3pt $\,p_{i\beta}\,$. Together with
the estimates on the exponential decay of the resonance function,
this will yield the sought estimate on $\,\im E\,$.
   \begin{proposition}
\label{lemlocfDthetarho}
Assume the conditions of theorem \ref{decay estimates thm}
and put $\,\rho_\star
:=\rho(\omega)\,$. Then there exists a number $\,c_{\eta}\,$
 such that \\
(i)
$
\|\,e^{-\rho}\,\im V^0_{i\beta}\,e^{-\rho}\,\|
\leq
c_{\eta}\, e^{-2\rho_\star}\,$,\vskip3pt
\noindent(ii)
$
\|\,\langle p \rangle ^{-1}e^{-\rho}\,\im W_{i\beta}e^{-\rho}
\langle p \rangle ^{-1}\|
\leq
dc_{\eta}\, e^{-2\rho_\star}
$ and \vskip3pt\noindent
(iii) there is a number $C_\eta$ such that
\label{form sense estimate for ImE}
$
\|\langle p \rangle ^{-1}
e^{-\rho}\,(\im H^j_{i\beta}-\,Z_{i\beta})\,e^{-\rho}
\langle p \rangle ^{-1}\|
\,\leq\, C_\eta e^{-2\rho_\star}\,.
$
  \end{proposition}
{\em Proof:} Since $\,\Sigma_{\alpha_0,\eta_0}\,$ and $\,\Sigma_{\beta,
\eta}\,$ are symmetric with respect to the real axis
we can choose the integration path $\,\partial\VV\,$ in
$\,\Sigma_{\alpha_0,\eta_0}\setminus \Sigma_{\beta,
\eta}\,$ invariant under complex conjugation. Using then
the Schwarz reflection principle, it is straightforward to compute
for a function $\,f\,$ obeying (a1)-(a2)
$$
\im f(D_{i\beta})
\,=\, {1\over 2i}\: {i\over 2\pi}\,
\int_{\partial\VV}\,
f(z)\, (r_{i\beta}(z)\!-\!r_{-i\beta}(z))\, dz\,.
$$
Now again by the norm convergence of the integral we have
$$
e^{-\rho}\,\im f(D_{i\beta})\,e^{-\rho}\,=\,
 {i\over 2\pi}\, \int_{\partial\VV}\,
f(z)\,e^{-\rho}\,\left({r_{i\beta}(z)\!-
\!r_{-i\beta}(z)\over 2i}\right)\,e^{-\rho}\, dz\,.
$$
Using equation (\ref{difference DTI}) with $\,\alpha =-\beta\,$ yields,
again omitting the argument $z$ for the resolvents, we get
\begin{eqnarray*}
\lefteqn{e^{-\rho}\left({r_{i\beta}\!-
\!r_{-i\beta}\over 2i}\right)e^{-\rho}}
\\ \\
&\, =\,  e^{-2\rho_\star}\left(\sin{\textstyle{\beta}}\,
z\,r_{-i\beta,-\rho}\chi_{\Omega_e}r _{i\beta,\rho}
+ 
\sin{\textstyle{\beta\over 2}} 
\left(e^{-i{\beta/2}}\chi_{\Omega_e}r _{i\beta,\rho}+
e^{i{\beta/2}}
r_{-i\beta,-\rho}\chi_{\Omega_e}\right)
\right)\,.
\end{eqnarray*}
Thus we obtain by Proposition~\ref{general reso p estimate}
   \begin{equation}\label{reso p difference}
\left\|e^{-\rho}\left({r_{i\beta}(z)\!-
\!r_{-i\beta}(z)\over 2i}\right)e^{-\rho}\right\|
\leq e^{-2\rho_\star}\left|\sin{\textstyle{\beta}}\right|\,
{ 2 \mbox{dist}(z,\Sigma_{\beta,\eta})+|z|\over
\mbox{dist}(z,\Sigma_{\beta,\eta})^2 }\,.
   \end{equation}
This justifies the statement for $\,V^0\,$ as in
Proposition~\ref{difference of the gammas}. \\
For (ii) we have
  \begin{eqnarray*}
\im \Bigl(p(b\!-\!1)p\Bigr)_{i\beta}
&\!=\!&
\,2\re\Bigl(
\im (p_{i\beta})\,(b\!-\!1)_{i\beta}p_{i\beta}\Bigl)
+
p_{-i\beta}\,\im (b\!-\!1)_{i\beta}\,p_{i\beta}
\\ \\
&\!=\!&
p\left(2\re\left(
{\im (p_{i\beta})\over p}\,(b\!-\!1)_{i\beta}{p_{i\beta}\over p}\right)
+
{p_{-i\beta}\over p}\,\im (b\!-\!1)_{i\beta}\,{p_{i\beta}\over p} \right)p
\,,
\end{eqnarray*}
yielding
 \begin{eqnarray*}
\Bigl\|\langle p \rangle^{-1} e^{-\rho}\,\im \Bigl(p
(b\!-\!1)p\Bigr)_{i\beta}\,e^{-\rho}\langle p \rangle^{-1}\Bigr\|
 &\!\!\le\!&\!
e^{-2\rho_\star}\,2\left\|(b\!-\!1)_{i\beta,\rho}\left\| +
\left\|e^{-\rho}\,\im b_{i\beta}\,e^{-\rho}
\right\|\right.\right.
\\ \\
 &\!\le\!&
c' _{\eta}\,d\,e^{-2\rho_\star}  \,,
   \end{eqnarray*}
for some number $\,c' _{\eta}\,$.
We have used here the fact that the imaginary part of $\,p_{i\beta}\,$ is
zero on $\,\Omega_i\,$ and  that $\,b\!-\!1= uf\,$, with $\,f\,$ obeying
(a1)--(a2) uniformly for $\,u\in [0,d],\ d\in (0,d_0)\,$. 
Thus we can apply the calculation used in (i) above.
This is also possible for $\,\im (V\!-\!V^0)_{i\beta}\,$.\\
Now (iii) is easy, since
$\,\im H^j_{i\beta}= \im (p^2_{i\beta}+ \JJ_j
W_{i\beta}\JJ^*_j +V^0_{i\beta})\,$.
Evidently we have
$$
\left\|\langle p \rangle^{-1} e^{-\rho}\,\im p^2_{i\beta}\,e^{-\rho}
\langle p \rangle^{-1} \right\| \,\leq\, e^{-2\rho_\star}
$$
The term
 $\,e^{-\rho}Z_{i\beta}e^{-\rho}$ is handled by
noting that $\,\im \widehat{H}^j _{i\beta} =
Q_j\im \left(p^2_{i\beta}+W _{i\beta}+
V^0_{i\beta} \right)Q_j\,$ and
that $\,\|\langle p \rangle \widehat{R}^j_{i\beta,\rho}
{W}_{i\beta\rho}\langle p \rangle^{-1}\|\leq 1\,$ by
Propositions~\ref{function B}(ii):
\begin{eqnarray*}\label{sRW}
\left\|
\langle p \rangle^{-1}e^{-\rho}Z_{i\beta}e^{-\rho}\langle p \rangle^{-1}
\right\|\!
&\!\!\leq\!\!&\!
2\left\|\langle p \rangle^{-1}e^{-\rho}\im W_{i\beta}e^{-\rho}
\langle p \rangle^{-1}\right\|
\!+\!
\left\|\langle p \rangle^{-1}e^{-\rho}\im H_{i\beta}e^{-\rho}
\langle p \rangle^{-1}\right\|
\\\nonumber \\
&\!\!\leq\!\!&\! (3 c_\eta d + c_\eta +1) e^{-2\rho_\star}\,. \quad\QED
 \end{eqnarray*}
Returning  to $\,\im E\,$ we know by general arguments that it cannot
be positive --- \cf\ \cite[Sec.XII.6]{RS},
so equation (\ref{fundamental eq for ImE}),
the above estimate, and Theorem~\ref{decay estimates thm}
yield
\begin{equation}
\label{bound on ImE}
0\,\leq\, -\im E\,\leq\,C_{\eta}\,e^{-2\rho_\star}
\left\{\|p\phi^j_{i\beta,\rho}\|^2 +
\tau\|\phi^j_{i\beta,\rho}\|^2\right\}
\,\leq\, {1\over 2}C_{\eta}\,(p_\star^2\!+\!\tau)
\,e^{-2\rho_\star}\,.
\end{equation}
The assertion of Theorem~\ref{full width
bound} now follows from the observation that 
$\,\tau\,$ and $\,p_\star\,$  are bounded
as $\,d\,$ tends to zero and that
$$
\exp\{-2\rho_\star\}\,=\,\exp\left\{-{2\pi\eta\over d}
\sqrt{2j-1}\,\left(1+\OO(\xi d)\right)\right\}\,.
$$

\setcounter{equation}{0}
\section{Proof of Theorem~2.3}
The proof uses the same ideas as the proof of Theorem~2.2 except that
due to the strengthened assumptions on the function $\,\gamma\,$,
we can allow now a boost function
$\,\rho\,$ with $\,\|\rho'\|_\infty\,$ exploiting asymptotically
the full width of the analyticity strip,
\ie $\,\|\rho'\|_\infty\,$ tending to $\,\eta_p\,$
as $\,d\,$ approaches zero.

The key to this is the representation of $\,f(D_{i\beta})\,$ below
when $\,f\,$ is a meromorphic function.
For the sake of simplicity, assume that $\,f\,$ has a single pair
of complex conjugated poles 
in $\,\Sigma_{\alpha_0,\eta_1} \setminus
\Sigma_{\alpha_0,\eta_0}\;$; an extension to any finite number is
straightforward. Let the order of these poles be $\,N\,$; for 
the proof of Theorem~2.3
we will have to consider several meromorphic functions made out of 
$\,\gamma\,$ with poles varying in order, not necessarily equal to 
$\,m\,$.
Without loss of generality, we also may suppose
that the poles lie on the imaginary axis at $\,z_p=i\eta_p\,$
and $\,{\overline z_p}\,$.
In view of the Schwarz reflection principle, it is sufficient to discuss
the behaviour of $\,f\,$ around the pole
in the upper half-plane and to translate the results by mirror
transformation to its counterpart; in particular, the integration contour
$\,\partial\VV\,$ in
the Dunford-Taylor integrals will always supposed to be symmetric
with respect to the real axis, \ie\ of the form
$\,\partial\VV:=\KK\cup\overline{\KK}\,$ with
a suitable upper branch $\,\KK\,$. By assumption, $\,f\,$ can be
expanded into its singular and regular part in a pierced neighbourhood
of $\,z_p\,$,
$$
f(p)\,=\, \sum_{k=1}^{N}{f_{-k}\over (p\!-\!z_p)^k}\,+\,
f_{reg}(p)
\,, \qquad 0<|p\!-\!z_p|<\eps\,.
$$
for some $\,\eps>0\,$. Let
$\,\KK\,$ now be passing above the pole  $\,z_p\,$, but lying entirely
inside $\,\Sigma_{\alpha_0,\eta_1}\,$ Then the residue theorem yields
the following
\begin{proposition}\label{prop: polar rep DTIO}
Let $\,f\,$ obey the same requirements as $\gamma$ in (a1)--(a2) and let
$\,f\,$ and
$\,\partial\VV\,$ be as above. Then
$$
f( D_{i\beta} )= \sum_{k=1}^{N}\left(
f_{-k}\,( D_{i\beta}\!-\!z_p)^{-k}+
\overline{f}_{-k}\,(D_{i\beta}\!-\!\overline{z}_p)^{-k}\right)
+ {i\over 2\pi} \int_{\partial\VV} f(z)( D_{i\beta}\!-\!z )^{-1}dz\,.
$$
\end{proposition}
This proposition together with Proposition~\ref{general reso p estimate}
yields immediately a bound on the boosted operator:
\begin{equation}\label{eq:7.1}
\bigl\|f( D_{i\beta,\rho} )\bigr\|\,\leq\,
{C\over (\eta_p -\|\rho'\|_\infty)^N}\,;
\end{equation}
note that the integral part can be uniformly bounded, since
the integration path $\,\KK\,$ can be kept at a finite distance,
independent of $\,d\,$, from the
horizontal $\,z=i\|\rho'\|_\infty\,$, so that this formula holds for all
$\,\|\rho'\|_\infty<\eta_p\,$ with an appropriate constant $\,C$.
In view of the basic decomposition (\ref{perturbation}) we thus have to
investigate how we can apply this formula to $\,b\!-\!1,\: V\!-\!V^0\,$,
and $\,V^0\,$ and how this conditions the maximal $\,\|\rho'\|_\infty\,$
to be chosen.
For this recall that we can interpret $\,b\,$ and $\,b\!-\!1\,$
as a simple rational function of $\,u\gamma\,$.
Choosing $\,\|\rho'\|_\infty=\eta_p\!-\!d^{1/(m+1)}\,$
implies $\,\|u\gamma_{i\beta,\rho} \|=\OO(d^{1/(m+1)})\,$
and $\left\|V^0_{i\beta,\rho}\right\| =\OO( d^{-2m/m+1})$,
the order of the pole of $\gamma$ being $\,m\,$.
Since $\,b\!-\!1= -u\gamma(2\!+\!u\gamma)b\,$ and
$\,V\!-\!V^0= V^0(b-1)
+ {1\over 2}u\gamma''b^{3/2}
- {5\over 4}u^2\gamma'^2b^2\,$, we obtain the bounds
applying the
above proposition and inequality (\ref{eq:7.1}) to the various
powers and powers of derivatives of
$\,\gamma\,$ observing that $\,\|b_{i\beta,\rho}\|\,$ is uniformly
bounded.

But before stating all the necessary bounds in a proposition
let us be more precise about the choice of
$\,\rho\,$. It shall be defined by formula (\ref{eq:def of rho}) with
$\,\eta\,$ replaced by
$\,\eta_p\!-\!d^{1/(m+1)}\,$ and $\,p_\star\,$ by
$\,p_\star'd^{-m/m+1}\,$ where  $\,p_\star' \,$ is a quantity uniformly
bounded with respect to $\,d\,$ to be fixed later. Note that
$\,\rho_\star\,$ still denotes $\,\rho(\omega)\,$. We also have to be
precise concerning the
weight $\,\langle p \rangle\,$:
$$
\langle p \rangle^2\,:=\, p^2+\tau,\quad \tau\,:=\, c_\gamma d^{-2m/m+1}\,.
$$
Recall that $\,\tau\,$ had been chosen to be a uniform bound on
$\left\|V^0_{i\beta,\rho}\right\|$.
As is confirmed in the next proposition this is again the case. 
All previous estimates involving $\,\tau\,$ used only this property
and remain thus valid. 
Notice that $\left\|V^0_{i\beta,\rho}\right\|$
does not depend on $\,p_\star' \,$, \cf (\ref{eq:7.1}).
\begin{proposition}
With the definitions above
and for $\,d\,$ small enough\\
(i) there exist numbers $c_\gamma$ and $c_b$ such that
$$
\sup_{0\leq -\beta\leq\beta_1}\left\|V^0_{i\beta,\rho}\right\| \leq c_\gamma
d^{{-2m\over m+1}}\,
\quad\mbox{ and }\quad
\bigl\|\langle p \rangle^{-1}W_{i\beta,\rho}\langle p \rangle^{-1}
\bigr\|
\le
{c_b}{d}^{{1\over m+1}}
\,.
$$
(ii)
There exists a constant $c$
such that
$$
\left\|\,e^{-\rho}\,\im V^0_{i\beta}\,e^{-\rho}\,
\right\|\leq
c\,{d}^{-{2m+1\over m+1}} \, e^{-2\rho_\star}\,\quad\mbox{ and }\quad
\,
\bigl\|\langle p \rangle^{-1}\,e^{-\rho}\,
\im W_{i\beta}\,e^{-\rho}\,\langle p \rangle^{-1}
\bigr\|\leq c\, e^{-2\rho_\star}
\,.
$$
\end{proposition}
{\em Proof:\/}
(i) The first statement is clear, the second
statement is obtained as in Proposition~\ref{W_thetarho estimate}.
We have here
$$
c_b =\max_{0\leq d\leq d_0}
\left\{ {d}^{-{1\over m+1}}\left(
\|(b\!-\!1)_{i\beta,\rho}\| +
 {\tau}^{-1}\|(V\!-\!V^0)_{i\beta,\rho}\|\right)\right\}\,, 
{\tau}^{-1}={{d}^{{2m\over m+1}}\over c_\gamma}\,.
$$
By the above discussion it is easy to see that  $\,c_b\,$ is uniformly
bounded in $\,d\,$. 
For (ii) we can use the same algebra as in
Proposition~\ref{lemlocfDthetarho}(i)--(ii).
It remains only to prove the proper localization of
the residual parts. We have,
switching back to $\,f\,$ as in Proposition~\ref{prop: polar rep DTIO}
and using the 
notation of Proposition~\ref{difference of the gammas},
for some $\, 1\leq k\leq N\,$
\begin{eqnarray*}
\lefteqn{\hskip-.7cm \im \left(
f_{-k}\,r_{i\beta}(z_p)^k+
\overline{f_{-k}}\,r_{i\beta}(\overline{z_p})^k\right)\,=\,}
&&\\ \\
&&\,=\,
\re f_{-k}\im \Bigl(r_{i\beta}(z_p)^k+{r_{-i\beta}({z_p})^k}^*\Bigr)+
\im f_{-k}\re \Bigl(r_{i\beta}(z_p)^k-{r_{-i\beta}({z_p})^k}^*\Bigr)
\\ \\
&&\,=\,
\re f_{-k}\im\Bigl(r_{i\beta}(z_p)^k- r_{-i\beta}(z_p)^k\Bigr)
+
\im f_{-k}\re\Bigl(r_{i\beta}(z_p)^k-r_{-i\beta}(z_p)^k\Bigr).
\end{eqnarray*}
The trivial identity 
$\,A^k\!-\!B^k=(A\!-\!B)A^{k-1}+B(A^{k-1}-B^{k-1})\,$ implies
$\,e^{-\rho}(A^k\!-\!B^k)e^{-\rho}=\sum_{\ell=0}^{k-1}
B^\ell_\rho e^{-\rho}(A\!-\!B)e^{-\rho}A^{k-1-\ell}_\rho$.
We obtain by Proposition~\ref{general reso p estimate}
\begin{eqnarray*}\left\|
e^{-\rho}\Bigl(r_{i\beta}(z_p)^k\!-\! r_{-i\beta}(z_p)^k\Bigr)
e^{-\rho}
\right\|
&\leq & k\,d^{-{k-1\over m+1}}
\left\|e^{-\rho}\Bigl(r_{i\beta}(z_p)\!-\! r_{-i\beta}(z_p)\Bigr)
e^{-\rho}\right\|\\
&\leq &
c'\,N\,d^{-{N+1\over m+1}}e^{-2\rho_\star};
\end{eqnarray*}
for the second inequality use (\ref{reso p difference}) and majorize
 $\,k\,$ by $N$. We explicitly have  $\,c':= |\sin\beta|(\eta_p\!+\!2
d^{1/(m+1)})\,$. 
Taking the appropriate $\,N\,$  for each of the functions concerned
yields the result.
\quad\QED\vskip4pt

Recall that the crucial equations to be justified are
(\ref{fundamental eq for decay estimate}) and
(\ref{fundamental eq for ImE}) which means the justification of the
existence of $\,B^j _{i\beta,\rho} (E)\,$ and thus of
$\widehat{R}^j _{i\beta,\rho} (E)\,$. Of course we still would like to
use the resolvent equation (\ref{reduced reso eq})
of Proposition~\ref{function B},
so we need
$\,
\bigl\|\langle p \rangle^{-1}\widehat{W}_{i\beta,\rho}
\widehat{R}^{0,j}_{i\beta,\rho}(E)\langle p \rangle\bigl\| \,<\,1$,
which is impossible
unless we make $\,\omega\,$ smaller.
We modify (\ref{dilation threshold}) by choosing
      \begin{equation} \label{modified dilation threshold}
\omega\,:=\, {\pi\over d}\, \sqrt{(2j\!-\!1)\left(1\!-\!\xi\,
{d}^{1/(m+1)}\right)}\,.
      \end{equation}
We follow the proof of Lemma~\ref{resolvent lemma} in Appendix~B.
Applying Proposition~\ref{double naught resolvent estimate}
with $\,\kappa=\xi\,{d}^{1/(m+1)}\,$ we get:
$$
\forall\,k\neq j\quad\, \|V^0_{i\beta,\rho} R^{00,k}_{i\beta}(z)\|
\,\leq\,
{c^{(1)}\tau\over 2|\sin{\beta}|}{d^2\over\xi{d}^{1/(m+1)} }
\,=\,
{c^{(1)}c_\gamma\over 2|\sin{\beta}|}{{d}^{1/(m+1)}\over\xi }\,.
$$
Thus for
$\,d\,$ small enough choosing $\,1\geq\xi\,{d}^{1/(m+1)}\geq 
c^{(1)}c_\gamma |\sin\beta|^{-1}{d}^{2/(m+1)}\,$ 
implies\\ $\,\|V^0_{i\beta,\rho} R^{00,k}_{i\beta}(z)\|\leq\,{1/2}\,$
and we obtain by 
the resolvent identity (\ref{miracle reso identity}) and (\ref{eq:
final reso esti})
\begin{equation}\label{eq: small R0j}
\left\|\,\langle p \rangle\widehat{R}^{0,j}_{i\beta,\rho}(E)
\langle p \rangle\right\|\,\le\,
{c^{(1)}\over |\sin\beta|\xi}\, {d}^{-1/(m+1)}
\,.
\end{equation}
Proceeding as in Proposition~\ref{function B} we need that 
$\,
\bigl\|\langle p \rangle^{-1}\widehat{W}_{i\beta,\rho}
\widehat{R}^{0,j}_{i\beta,\rho}(E)\langle p \rangle\bigl\| \,\leq\,1/2$,
which is the case if 
$\,\xi\geq {2c^{(1)}c_b |\sin\beta|^{-1}}\,$. Consequently 
$$
\left\|\,\langle p \rangle\widehat{R}^{j}_{i\beta,\rho}\langle p \rangle
\,\right\|\,\le\,
{2c^{(1)}\over|\sin\beta|\xi}\, {d}^{-1/(m+1)}
\quad\mbox{ and }\quad
\|\langle p \rangle^{-1}B^j_{i\beta,\rho}(E)\langle p \rangle^{-1}\|
\,\le\,
c_b\, {d}^{1/(m+1)}\,.
$$
The
inequality of Proposition~\ref{form sense estimate for Agmon} can be
formulated now as
$$
\re\left(\,H^j_{i\beta,\rho}
 -B^j_{i\beta,\rho}(E)-E\right)
\,\geq\,
\left(\cos 2\beta\!-\!2c_b\, {d}^{1/(m+1)}\right)
\left(p^2\! -\! {{p'_\star}^2\over 2}{d}^{-2m/(m+1)} \right)
$$
with the new, but nevertheless uniformly bounded (in $\,d\,$)
$$
{{p'_\star}^2\over 2}\,:=\,{c_\gamma \!+\!2c_b\,c_\gamma{d}^{1/(m+1)}\! +\!
{d}^{2m/(m+1)}|\re E\!-\!E_j|
\over \cos 2\beta\!-\!2c_b\, {d}^{1/(m+1)}}\,;
$$
the inequality is, of course, to be interpreted
in the form sense on $\,\DD(p^2\otimes I)\,$.
Thus for all sufficiently small $\,d\,$
the formula (\ref{eq:reso fct estimate}) remains valid
with $\,p_\star=p'_\star{d}^{-m/(m+1)} \,$
when $\,\eta\,$ is changed to
$\,\eta_p\!-\!{d}^{1/(m+1)}\,$ in the definition
(\ref{eq:def of rho}) of $\rho$.

Since also the algebra used for Proposition~\ref{form
sense estimate for ImE}(iii) can be applied without change, we
just need to substitute  corresponding constants to arrive
at the following inequality replacing (\ref{bound on ImE})
$$
0\,\geq\, \im E\,\geq\, - C\, ({p'}^2_\star +c_\gamma)\,d^{-4+3/{(m+1)}}
\,e^{-2\rho_\star}\,,
$$
for some constant $\,C\,$.
To conclude the proof, it remains to expand $\,\rho_\star\,$:
   \begin{eqnarray*}
\rho_\star &\!=\!& (\eta_p\!-\! {d}^{1/(m+1)}) \left({\pi\over d}\,
\sqrt{(2j\!-\!1)\left(1\!-\!\xi\,{d}^{1/(m+1)}\right)}\,
 -{p'_\star}{d}^{-m/(m+1)}\right) \\ \\
&\!=\!& {\pi\eta_p\over d}\, \sqrt{2j\!-\!1}\,
\left(1\!+\!\OO({d}^{1/(m+1)})\,\right)\,,
   \end{eqnarray*}
and to notice that negative powers of $d$ in the prefactor
have no significance and
can be absorbed in the error term of the exponential decay rate.
\appendix
\renewcommand{\theequation}{\Alph{section}.\arabic{equation}}
\setcounter{equation}{0}
\section{Proof of Lemma~3.4}
{\em Proof of Lemma~\ref{DTI properties}:\/}
(i) By hypothesis $\,(T-z)^{-1}$ is bounded on the integration path and
$\,f\,$ decays rapidly enough to make the integral
converge in operator norm.
Furthermore, the integral does not depend on the path, since both the
resolvent of $\,T\,$ as a function of $\,z\,$ and $\,f\,$ are analytic in the
considered region. (ii) Since $\,f(T)\,$ is bounded, it suffices to show that
$\,(f(T)u,v)= (f_{sp}(T)u,v)\,$ holds for all $\,u,v\in L^2(\R)\,$,
where $\,f_{sp}(T)\,$ denotes the operator
defined by the spectral theorem. One has
   \begin{eqnarray*}
(f_{sp}(T)u,v) &\!=\!& \int_{\R} f(\lambda)\: d(E_\lambda u,v) \\ \\
&\!=\!& \int_{\R} d(E_\lambda u,v)
{i\over 2\pi} \int_{\partial \VV}{ f(z)\over
\lambda \!-\!z}\:dz\, \\ \\
&\!=\!&{i\over 2\pi}
\int_{\partial \VV} dz
\int_{\R}{ f(z)\over \lambda \!-\!z}\: d(E_\lambda u,v)\,\\ \\
&\!=\!& {i\over 2\pi}
\int_{\partial \VV} f(z)((T\!-\!z)^{-1}u,v)\,dz \,=\,(f(T)u,v)\,,
   \end{eqnarray*}
where in the third and the last step we have employed the Fubini theorem. \\
(iii) It follows from the norm convergence of the integral that
$$
B\, {i\over 2\pi}\, \int_{\partial \VV}\,
f(z) (T\!-\!z)^{-1}\, dz\,B^{-1}\,=\,
{i\over 2\pi}\, \int_{\partial \VV}\,
f(z)\, B(T\!-\!z)^{-1}B^{-1}\, dz\,=\,
f(BTB^{-1})\,.
$$
(iv)
The operators $\,U_{\theta}\,$ are unitary for $\,\theta
\in \R\,$, and therefore $\,\Big(f(T)\Big)_{\theta}= f(T_{\theta})\,$ 
by  (iii).
Furthermore, the resolvent
$\,(T_{\theta}\!-\!z)^{-1}\,$ is by hypothesis uniformly bounded
on $\,\partial \VV\,$ for all $\,{\theta}\in \SS_\alpha$,
and  analytic in $\,{\theta}\,$.
Thus the analyticity follows by the convergence in operator norm of
the integral, since the limit function of a uniformly
convergent sequence of analytic functions is analytic (see \eg
\cite[Thm. 9.12.1]{Di}).
\quad \QED
\vspace{3mm}

\setcounter{equation}{0}
\section{Proof of Lemma 4.2}
We need the following
        \begin{proposition}\label{double naught resolvent estimate}
Let $\,R^{00,k}_{i\beta}(z):= (p_{i\beta}^2 +E_k - z)^{-1}$,
where $\,0<|\beta|\leq \min\{{\pi\over 4},\alpha\}\,$.
Let $\omega\geq 0$ be defined by
$\displaystyle
\,\omega^2 =(1\!-\!\kappa)(E_j-E_{j-1})\,
$,
where $\,\kappa\in (0,{1}]\,$ and  $\,c^{(1)}:=8\sqrt{3}\, $.
Then
for all  $\,(\kappa,d)\,$  such that
\begin{equation}\label{first domain kappa d}
{1}\geq\,\kappa\,\geq\,{\tau\over \pi^2|\sin\beta|}\,d^2\,,
\end{equation}
all $z\in \DD_\Gamma\,$ and $\,\ell=0,1,2\,$ one has
    \begin{description}
\item{(i)}
$\displaystyle
\quad\|\langle p \rangle^{\ell}
R_{i\beta}^{00,j-1}(z)\|\le {c^{(1)}\over 2|\sin\beta|\kappa}\,
d^{2-\ell}\,$ and
\item{(ii)}
$\displaystyle
\quad \|\langle p \rangle^{\ell}
R_{i\beta}^{00,k}(z)\|\le {c^{(1)}\over 2|\sin\beta|}\,d^{2-\ell}\quad
\forall\, k\ne j,\, j\!-\!1\, $.
   \end{description}
        \end{proposition}
{\em Proof 
:\/}
We first estimate $ R^{00,k}_{i\beta}(z),\;k\neq j\,$.
 Define $\,\zeta:=z-E_{j}\,$
and $\,\Delta_{j,k}:=E_j-E_k\,$.
If
\begin{equation}\label{first domain zeta d}
-{1\over 2}\kappa\Delta_{j,j-1}\leq \re\zeta\leq 0
\quad{\rm and}\quad
\im e^{-i\beta}\zeta \geq {1\over 2}\kappa\Delta_{j,j-1}\sin\beta,
\end{equation}
then one obtains by simple geometric considerations
$$
\left\|R^{00,k}_{i\beta}(z)\right\|
\,\leq\,{2\over |\sin\beta|\Delta_{j,k}},\,k<j-1,\quad\quad
\left\|R^{00,k}_{i\beta}(z)\right\|
\,\leq\,{1\over \Delta_{k,j}},\,k>j,
$$
and
$$
\left\|R^{00,j-1}_{i\beta}(z)\right\|
\,\leq {2\over |\sin\beta|\kappa \Delta_{j,j-1}}.
$$
The condition
$
\kappa\geq (\pi^2|\sin\beta|)^{-1}{\|V^0\|}d^2
$
is sufficient to ensure that
$\DD_\Gamma$ is for all $j\geq 2$
contained in the domain described
by (\ref{first domain zeta d}) and we have, of course,
$\,\|V^0\|\leq\tau\,$. Thus the case $\,\ell=0\,$ is proven
noticing that $\,|\Delta_{k,j}|^{-1},k\neq j\,$, is uniformly bounded
by $\,{(3\pi^2)^{-1}d^2}\,$.

To treat the case $\,\ell=2\,$
we write
$$\,\langle p \rangle R^{00,k}_{i\beta}(z)\langle p \rangle
={(p^2+\tau)\over (p_{i\beta}^2+\tau)}
\left(1+(\tau-E_k+z)R^{00,k}_{i\beta}(z)
\right)\,.
$$
The first factor is uniformly bounded by
$\sqrt{3}$ for $\,|\beta|\leq {\pi\over 4}\,$.
Again simple geometric considerations suffice to bound  
the term $\,|E_k-z|\|R^{00,k}_{i\beta}(z)\|\,$. Note that
the overall constant is  made $d$ independent by condition
(\ref{first domain kappa d}).

The remaining case $\,\ell=1\,$ is handled by
the inequality
$$
\left\|\langle p \rangle R^{00,k}_{i\beta}(z)\right\|
\,\leq \,\left\|R^{00,k}_{i\beta}(z)\right\|^{1/2}
\left\|\langle p \rangle^2 R^{00,k}_{i\beta}(z)\right\|^{1/2}\,.
 \quad \QED
$$
\vskip3pt\noindent
{\em Proof of Lemma~\ref{resolvent lemma}:\/}
For $\ell=0$ observe that  by replacing $\kappa$ with $\xi d$
one has
\begin{equation}\label{eq:aux}
\, \|V^0_{i\beta,\rho} R^{00,k}_{i\beta}(z)\|
\,\leq\,
\|V^0_{i\beta,\rho}\|\,\| R^{00,k}_{i\beta}(z)\|
\,\leq\,
{c^{(1)}\tau\over 2|\sin{\beta}|}{d\over\xi }
\,\leq\,
{1\over 2}
\end{equation}
by the condition on $\xi$,
uniformly for all $\,z\in \DD_{\Gamma}\,$ and $\,k\ne j$,
so that by the above estimates the bounds
follow immediately in this case.
To prove the estimates in the case $\ell=1$ we use the
resolvent identity
\begin{equation}\label{miracle reso identity}
\langle p \rangle R^{0,k}_{i\beta,\rho}\langle p \rangle \,=\,
\langle p \rangle R^{00,k}_{i\beta}\langle p \rangle +
 \langle p \rangle R^{00,k}_{i\beta}V^0_{i\beta,\rho} \left( 1\!
+\! R^{00,k}_{i\beta}V^0_{i\beta,\rho}\right)^{-1}
R^{00,k}_{i\beta}\langle p \rangle \,.
\end{equation}
This yields
\begin{eqnarray}\nonumber
\left\|\langle p \rangle R^{0,k}_{i\beta,\rho}(z)\langle p \rangle \right\|
&\!\leq\!&
\left\|\langle p \rangle^2 R^{00,k}_{i\beta}(z)\right\|
+\left\|\langle p \rangle R^{00,k}_{i\beta}(z)\right\|^2
{\|V^0_{i\beta,\rho}\|\over 1 - \|V^0_{i\beta,\rho}\|\,
\| R^{00,k}_{i\beta}(z)\|}\\\nonumber \\\nonumber 
&\!\leq\!&
\left\|\langle p \rangle^2 R^{00,k}_{i\beta}(z)\right\|
+\left\|\langle p \rangle^2 R^{00,k}_{i\beta}(z)\right\|
{\|V^0_{i\beta,\rho}\|\,
\| R^{00,k}_{i\beta}(z)\|\over 1 - \|V^0_{i\beta,\rho}\|\,\|
 R^{00,k}_{i\beta}(z)\|}\\\nonumber  \\ \label{eq: final reso esti}
&\!\leq\!& 2\left\|\langle p \rangle^2 R^{00,k}_{i\beta,\rho}(z)\right\|
\end{eqnarray}
using in the last step (\ref{eq:aux}), and in the second to last
step the fact that $\, R^{00,k}_{i\beta}(z)\,$ is a multiplication
operator. So the bounds follow again easily.
\quad\QED
\vspace{8pt}\noindent

\noindent
{\bf Acknowledgment.} \quad The work has been done during the visits
of P.~D. and B.~M. to the Nuclear Physics Institute, AS CR, and of
P.E. to Centre de Physique Th\'eorique, Marseille--Luminy and
Universit\'e de Toulon et du Var; the authors express their gratitude
to the hosts. A partial support by the grants AS CR No. 148409 and
GA CR No. 202--93--1314 is  also gratefully acknowledged.

   \end{document}